\documentclass[final]{article}

\usepackage{amsfonts}
\usepackage{amsmath}
\usepackage{amssymb}
\usepackage{mathrsfs}

\usepackage{showkeys}
\usepackage{graphicx}
\usepackage{color}
\usepackage{afterpage}
\usepackage{verbatim}
\usepackage{afterpage}
\usepackage{amsthm}
\usepackage[normalem]{ulem}
\usepackage{bm}
\usepackage{hyperref}
\hypersetup{
   unicode=true,
   pdftoolbar=false,
   pdfmenubar=false,
   pdffitwindow=false,
   pdfnewwindow=true,
   pdfkeywords={},
   colorlinks=true,
   linkcolor=black,
   citecolor=black,
   filecolor=black,
   urlcolor=black
}

\usepackage{caption}
\usepackage{subcaption}
\usepackage{tikz}
\usepackage{pgfplotstable}
\usetikzlibrary{spy,calc}
\usepackage{hyperref}

\newif\ifblackandwhitecycle
\gdef\patternnumber{0}

\pgfkeys{/tikz/.cd,
    zoombox paths/.style={
        draw=orange,
        very thick
    },
    black and white/.is choice,
    black and white/.default=static,
    black and white/static/.style={ 
        draw=white,   
        zoombox paths/.append style={
            draw=white,
            postaction={
                draw=black,
                loosely dashed
            }
        }
    },
    black and white/static/.code={
        \gdef\patternnumber{1}
    },
    black and white/cycle/.code={
        \blackandwhitecycletrue
        \gdef\patternnumber{1}
    },
    black and white pattern/.is choice,
    black and white pattern/0/.style={},
    black and white pattern/1/.style={    
            draw=white,
            postaction={
                draw=black,
                dash pattern=on 2pt off 2pt
            }
    },
    black and white pattern/2/.style={    
            draw=white,
            postaction={
                draw=black,
                dash pattern=on 4pt off 4pt
            }
    },
    black and white pattern/3/.style={    
            draw=white,
            postaction={
                draw=black,
                dash pattern=on 4pt off 4pt on 1pt off 4pt
            }
    },
    black and white pattern/4/.style={    
            draw=white,
            postaction={
                draw=black,
                dash pattern=on 4pt off 2pt on 2 pt off 2pt on 2 pt off 2pt
            }
    },
    zoomboxarray inner gap/.initial=5pt,
    zoomboxarray columns/.initial=2,
    zoomboxarray rows/.initial=2,
    subfigurename/.initial={},
    figurename/.initial={zoombox},
    zoomboxarray/.style={
        execute at begin picture={
            \begin{scope}[
                spy using outlines={%
                    zoombox paths,
                    width=\imagewidth / \pgfkeysvalueof{/tikz/zoomboxarray columns} - (\pgfkeysvalueof{/tikz/zoomboxarray columns} - 1) / \pgfkeysvalueof{/tikz/zoomboxarray columns} * \pgfkeysvalueof{/tikz/zoomboxarray inner gap} -\pgflinewidth,
                    height=\imageheight / \pgfkeysvalueof{/tikz/zoomboxarray rows} - (\pgfkeysvalueof{/tikz/zoomboxarray rows} - 1) / \pgfkeysvalueof{/tikz/zoomboxarray rows} * \pgfkeysvalueof{/tikz/zoomboxarray inner gap}-\pgflinewidth,
                    magnification=3,
                    every spy on node/.style={
                        zoombox paths
                    },
                    every spy in node/.style={
                        zoombox paths
                    }
                }
            ]
        },
        execute at end picture={
            \end{scope}
            \node at (image.north) [anchor=north,inner sep=0pt] {\subcaptionbox{\label{\pgfkeysvalueof{/tikz/figurename}-image}}{\phantomimage}};
            \node at (zoomboxes container.north) [anchor=north,inner sep=0pt] {\subcaptionbox{\label{\pgfkeysvalueof{/tikz/figurename}-zoom}}{\phantomimage}};
     \gdef\patternnumber{0}
        },
        spymargin/.initial=0.5em,
        zoomboxes xshift/.initial=1,
        zoomboxes right/.code=\pgfkeys{/tikz/zoomboxes xshift=1},
        zoomboxes left/.code=\pgfkeys{/tikz/zoomboxes xshift=-1},
        zoomboxes yshift/.initial=0,
        zoomboxes above/.code={
            \pgfkeys{/tikz/zoomboxes yshift=1},
            \pgfkeys{/tikz/zoomboxes xshift=0}
        },
        zoomboxes below/.code={
            \pgfkeys{/tikz/zoomboxes yshift=-1},
            \pgfkeys{/tikz/zoomboxes xshift=0}
        },
        caption margin/.initial=4ex,
    },
    adjust caption spacing/.code={},
    image container/.style={
        inner sep=0pt,
        at=(image.north),
        anchor=north,
        adjust caption spacing
    },
    zoomboxes container/.style={
        inner sep=0pt,
        at=(image.north),
        anchor=north,
        name=zoomboxes container,
        xshift=\pgfkeysvalueof{/tikz/zoomboxes xshift}*(\imagewidth+\pgfkeysvalueof{/tikz/spymargin}),
        yshift=\pgfkeysvalueof{/tikz/zoomboxes yshift}*(\imageheight+\pgfkeysvalueof{/tikz/spymargin}+\pgfkeysvalueof{/tikz/caption margin}),
        adjust caption spacing
    },
    calculate dimensions/.code={
        \pgfpointdiff{\pgfpointanchor{image}{south west} }{\pgfpointanchor{image}{north east} }
        \pgfgetlastxy{\imagewidth}{\imageheight}
        \global\let\imagewidth=\imagewidth
        \global\let\imageheight=\imageheight
        \gdef\columncount{1}
        \gdef\rowcount{1}
        
    },
    image node/.style={
        inner sep=0pt,
        name=image,
        anchor=south west,
        append after command={
            [calculate dimensions]
            node [image container,subfigurename=\pgfkeysvalueof{/tikz/figurename}-image] {\phantomimage}
            node [zoomboxes container,subfigurename=\pgfkeysvalueof{/tikz/figurename}-zoom] {\phantomimage}
        }
    },
    color code/.style={
        zoombox paths/.append style={draw=#1}
    },
    connect zoomboxes/.style={
    spy connection path={\draw[draw=none,zoombox paths] (tikzspyonnode) -- (tikzspyinnode);}
    },
    help grid code/.code={
        \begin{scope}[
                x={(image.south east)},
                y={(image.north west)},
                font=\footnotesize,
                help lines,
                overlay
            ]
            \foreach \x in {0,1,...,9} { 
                \draw(\x/10,0) -- (\x/10,1);
                \node [anchor=north] at (\x/10,0) {0.\x};
            }
            \foreach \y in {0,1,...,9} {
                \draw(0,\y/10) -- (1,\y/10);                        \node [anchor=east] at (0,\y/10) {0.\y};
            }
        \end{scope}    
    },
    help grid/.style={
        append after command={
            [help grid code]
        }
    },
}

\newcommand\phantomimage{%
    \phantom{%
        \rule{\imagewidth}{\imageheight}%
    }%
}
\newcommand\zoombox[2][]{
    \begin{scope}[zoombox paths]
        \pgfmathsetmacro\xpos{
            (\columncount-1)*(\imagewidth / \pgfkeysvalueof{/tikz/zoomboxarray columns} + \pgfkeysvalueof{/tikz/zoomboxarray inner gap} / \pgfkeysvalueof{/tikz/zoomboxarray columns} ) + \pgflinewidth
        }
        \pgfmathsetmacro\ypos{
            (\rowcount-1)*( \imageheight / \pgfkeysvalueof{/tikz/zoomboxarray rows} + \pgfkeysvalueof{/tikz/zoomboxarray inner gap} / \pgfkeysvalueof{/tikz/zoomboxarray rows} ) + 0.5*\pgflinewidth
        }
        \edef\dospy{\noexpand\spy [
            #1,
            zoombox paths/.append style={
                black and white pattern=\patternnumber
            },
            every spy on node/.append style={#1},
            x=\imagewidth,
            y=\imageheight
        ] on (#2) in node [anchor=north west] at ($(zoomboxes container.north west)+(\xpos pt,-\ypos pt)$);}
        \dospy
        \pgfmathtruncatemacro\pgfmathresult{ifthenelse(\columncount==\pgfkeysvalueof{/tikz/zoomboxarray columns},\rowcount+1,\rowcount)}
        \global\let\rowcount=\pgfmathresult
        \pgfmathtruncatemacro\pgfmathresult{ifthenelse(\columncount==\pgfkeysvalueof{/tikz/zoomboxarray columns},1,\columncount+1)}
        \global\let\columncount=\pgfmathresult
        \ifblackandwhitecycle
            \pgfmathtruncatemacro{\newpatternnumber}{\patternnumber+1}
            \global\edef\patternnumber{\newpatternnumber}
        \fi
    \end{scope}
}

\numberwithin{equation}{section}

\newtheorem*{theo*}{Theorem}
\newtheorem*{corol*}{Corollary}

\newcommand{\specialcell}[2][l]{%
  \begin{tabular}[#1]{@{}l@{}}#2\end{tabular}}

\newcommand{\FT}{\mathcal{F}}
\newcommand{\bb}[1]{{\bm #1}}

\newcommand{\Radon}{\mathscr{R}}
\newcommand{\Back}{\mathscr{B}}

\newcommand{\R}{\mathbb{R}}

\newcommand{\ds}{\displaystyle}

\newcommand{\implica}{\Rightarrow}

\newcommand{\minimize}{\displaystyle \operatornamewithlimits{minimize}}

\graphicspath{  {figures/}}

\begin{document}

\title{Fast Backprojection Techniques for High Resolution Tomography}

\author{Nikolay Koshev,  Elias S. Helou, Eduardo X. Miqueles}

\maketitle

\begin{abstract}
	Fast image reconstruction techniques are becoming important with the increasing number of scientific cases 
	in high resolution micro and nano tomography. The processing of the large scale three-dimensional data
	demands new mathematical tools for the tomographic reconstruction task because of the big computational complexity of most current algorithms as the sizes of tomographic data grow with the development of more powerful acquisition hardware and
	more refined scientific needs. In the present paper we propose a new fast back-projection operator for the
	processing of tomographic data and compare it against other fast reconstruction techniques.
\end{abstract}

\section{Introduction}

Tomographic imaging is a very powerful instrument of non-destructive research and control of the internal structure of non-opaque objects. An important branch of tomographic techniques is \emph{transmission tomography}, which can be used at nano, micro and macro resolution levels.
For further consideration we describe in general the basic principles of transmission tomography from parallel rays, and define some notations.

\begin{figure}[t]
\centering
\includegraphics[width=\textwidth]{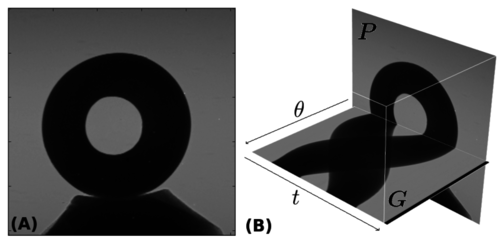}
\caption{(A) Projection (or \emph{frame}) for a cylindrical sample obtained with a \textsc{ccd} camera (B) 
Three-dimensional representation of the measured data: $P$ is the measured frame and $G$ is the sinogram image at a given row of the area detector.}
\label{fig:frame}
\end{figure}

Physically, all types of transmission tomography are based on registering the 
energy loss or/and intensity loss of the incoming electromagnetic wave (\textsc{x}-rays for instance), 
after passing through the object under investigation also referred to here as \emph{sample}). 
In our case, we consider that \textsc{x}-rays generated from a synchrotron light source hit the object under investigation 
determining a projection image (also referred as \emph{frame}) at a \textsc{ccd} (charge coupled device) camera. A typical 
dataset is shown in Figure \ref{fig:frame}.A, where a
high-resolution frame $P$ gathered using the \textsc{x}-rays source is shown, with dimensions $2048\times 2048$.
After half-rotation of the sample on the rotation axis, we obtain a cubic dataset as shown in Figure \ref{fig:frame}.B. 
Each slice of this dataset give us an image, which is called \emph{sinogram}, and that will be used as input to an appropriate inversion algorithm in order to reconstruct the slice of the sample.

We introduce the cartesian coordinate system in the plane of a given slice of the object. 
Let the function $f(\bb x) \in U$ be the \emph{feature function}, i.e., a 
function which depends on the internal structure of the object in the plane of the slice 
and which defines the linear absorption coefficient of the sample. Set $U$, referred to 
here as the \emph{feature space}, is a Schwartz space $\mathcal S(\R^2)$. 

\begin{figure}
\centering
\includegraphics[scale=0.35]{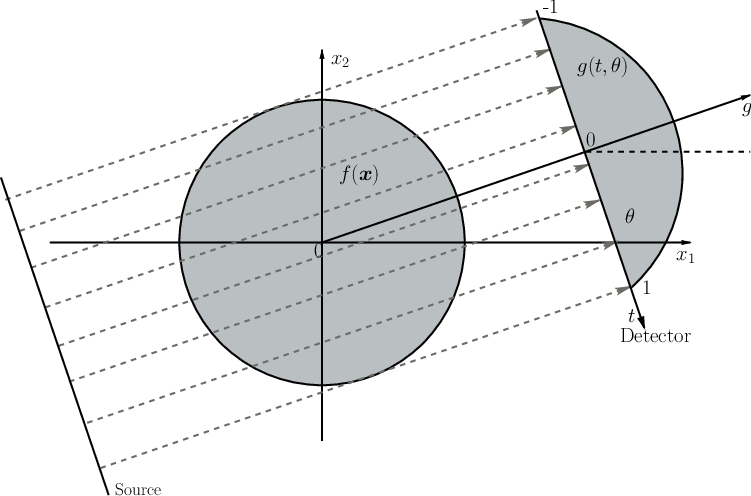}
\caption{Geometry of incoming x-rays for parallel tomography.}
\label{fig:scheme}
\end{figure}

A given frame (see Figure \ref{fig:frame}.A) represents the integral of $f(\bb x)$ over 
straight lines passing through the sample and perpendicular to 
the detector's plane. One row of each of these frame images contains the integrals relevant to a slice of the object,
and orthogonal to the rotation axis. Let us introduce an axis $t$ over the detector's row. It is clear that for each angle 
$\theta$ (see Figure \ref{fig:scheme}), such a row is mathematically determined by
\begin{equation}\label{P_theta}
	g(\theta, t) \equiv g_{\theta}(t) = \int\limits_{L(\theta, t)} f(\bb x) \mathrm d s = \int\limits_{\mathbb{R}^2} 
		f(\bb x) \delta( \bb x \cdot \bb \xi_\theta -t ) \mathrm d\bm x, 
\end{equation}
where $L(\theta, t)$ is a straight line defining the x-ray path,
\begin{equation}\label{L}
	L{(\theta, t)} = \big\{ \bb x \in \R^2: \bb x \cdot \bb \xi_\theta = t \big\} , \ \ \ \bb \xi_\theta = (\cos\theta,\sin\theta)^T.
\end{equation}
From (\ref{P_theta}) we have a linear operator acting on the feature function $f$, i.e.,
$\Radon \colon f \in U \mapsto g \in V$, which is called the \emph{Radon transform}. Space $V$
is the Schwartz space $\mathcal S(\R_+ \times [0,\pi])$. The operator 
$\Back \colon V \to U$ defined as
\begin{equation}\label{B}
	b(\bb x) = \Back g(\bb x) = \int_{[0,\pi]} g(\bb x\cdot \bb \xi_\theta, \theta) \mathrm d\theta, 
\end{equation} 
is defined as the \emph{backprojection} operator, and is the adjoint of $\Radon$ in the following
sense
\begin{equation}\label{adjoint}
	\int\limits_{\R_+ \times [0,\pi]} \Radon f(t, \theta) g(t, \theta) \mathrm d t \mathrm d\theta= 
	\int\limits_{\mathbb{R}^2} f(\bb x) \Back g(\bb x) \mathrm d \bb x,  
\end{equation} 
More about the theory of the integral operators $\{\Radon, \Back\}$ can be found on \cite{deans, kak_slaney, nattw, helgason}.

\begin{figure}
	\begin{center}
		\begin{tabular}{cc}
		(a) & (b) \\
		\begin{tikzpicture}
	        \begin{axis}[width=6cm, compat=newest,
	        ticks=none,enlargelimits=false, axis on top, scale only axis, 
	        axis equal image, xlabel={$|\bm x_1| \leq 1$}, ylabel={$|\bm x_2|\leq 1$}]
	        \addplot graphics 
	        [xmin=0,xmax=256,ymin=0,ymax=256]{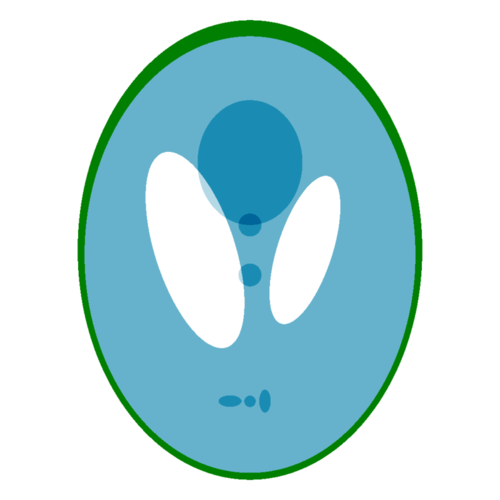};
	        \end{axis}
	        \end{tikzpicture} &
		\begin{tikzpicture}
	        \begin{axis}[width=6cm, compat=newest,
	        ticks=none,enlargelimits=false, axis on top, scale only axis, 
	        axis equal image, ylabel={$t\in [-\sqrt 2, \sqrt 2]$}, xlabel={$\theta\in [0,\pi]$ }]
	        \addplot graphics 
	        [xmin=0,xmax=256,ymin=0,ymax=256]{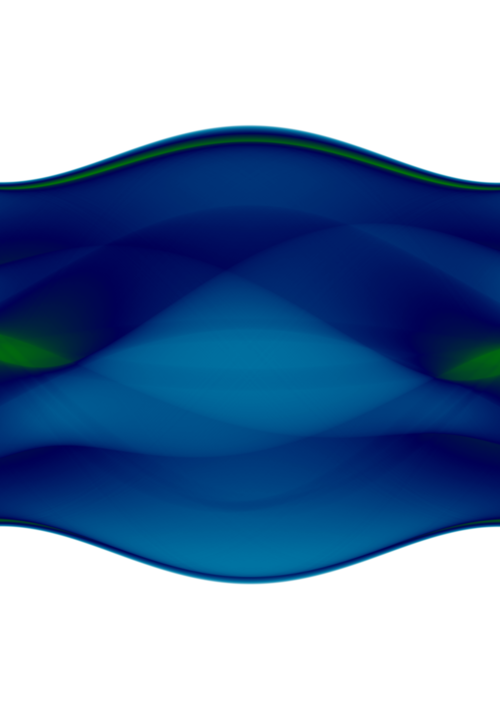};
	        \end{axis}
 	        \end{tikzpicture} \\
		(c) & (d) \\
		\begin{tikzpicture}
	        \begin{axis}[width=6cm, compat=newest,
	        ticks=none,enlargelimits=false, axis on top, scale only axis, 
	        axis equal image, ylabel={$t \in [0,\sqrt 2]$}, xlabel={$\theta \in [0,2\pi]$}]
	        \addplot graphics 
	        [xmin=0,xmax=256,ymin=0,ymax=256]{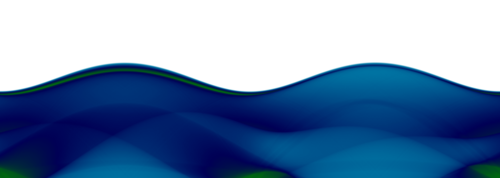};
	        \end{axis}
		\end{tikzpicture} &
		\begin{tikzpicture}
	        \begin{axis}[width=6cm, compat=newest,
	        ticks=none,enlargelimits=false, axis on top, scale only axis, 
	        axis equal image, ylabel={$\rho \in (-\infty, 0]$}, xlabel={$\theta \in [0,2\pi]$}]
	        \addplot graphics 
	        [xmin=0,xmax=256,ymin=0,ymax=256]{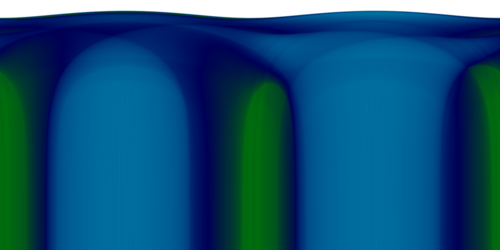};
	        \end{axis}
		\end{tikzpicture}
		\end{tabular}
	\end{center}
	\caption{\small (a) Shepp-Logan feature function $f(\bb x)$ and his associated sinograms, in different
	coordinate systems: (b) Semi-polar coordinates $[g]_{\sf s}$, (c) Polar coordinates $[g]_{\sf P}$ and (d) Log-polar coordinates $[g]_{\sf L}$.}
	\label{fig:phantom_sino}
\end{figure}

At this point, it is convenient to introduce some notations. We first introduce the notations for the representation of feature function $f: \R^2 \to \R$ in different coordinate systems, and their respective jacobians:
\begin{description}
	\item[(a)] \emph{Pr{\"u}fer coordinates (see \cite{prufer})}: $\bb x = p(\mu)\bb \xi_{\theta}$, $\mathrm d \bb x 
		= |p'(\mu)p(\mu)| \mathrm d\mu \mathrm d\theta$. The representation is denoted by 
		$[f]_{\mathsf{Pr}}(\mu, \theta)$. Function $p$ will always be well defined within the context by special notation as follows;
	\item[(b)] \emph{Log-polar coordinates}: particular case of Pr{\"u}fer coordinates when 
		$p(\mu) = e^{\mu}$. Here, $\mathrm d \bb x = e^{2\mu} \mathrm d\mu \mathrm d\theta$. 
		The representation is denoted by $[f]_{\mathsf L}(\mu, \theta)$;
	\item[(c)] \emph{Semi-polar coordinates}: particular case of Pr{\"u}fer coordinates when 
		$p(\mu) = \mu$. Here, $\mathrm d \bb x = \mu \mathrm d\mu \mathrm d\theta$. 
		The representation is denoted by $[f]_{\mathsf P}(\mu, \theta)$.
	\item[(d)] \emph{Sinogram} coordinates are similar to the semi-polar 
		coordinates and, in fact, can be obtained by flipping the angles 
		$\theta \in [\pi, 2\pi)$ to the negative part of the $t$-axis ,
		so that $t \in[-1, 1]$.
\end{description}
Using the above notation, function $g$ in (\ref{P_theta}) can be written as $[g]_{\sf P}(t,\theta)$ in order to indicate semi-polar coordinates. An example, using the well-known \emph{Shepp-Logan phantom} \cite{shepplogan} is presented on Figure \ref{fig:phantom_sino}. The sinogram of the Shepp-logan feature function $f$ is presented in the sinogram coordinate system mentioned above.

Our goal in the present paper is to present a fast method for the computation of the backprojection image $b\in U$ 
for a given sinogram $g \in V$. The computation of the tomographic image from sinogram data depends on the backprojection
operator $\Back$, which bears the major computational cost of reconstruction methods: $O( N^3 )$ for images with $N^2$
pixels and $N$ projections.


For a high-resolution tomographic synchrotron experiment, the amount of data at a micro-tomography setup is considerably  
large for today's computational standards, mainly because of this asymptotic floating point operations (flops) count. Indeed, at the 
Brazilian National Synchrotron Light Source (\textsc{lnls}) one wishes to obtain $2048$ reconstructions images with $2048\times 2048$ 
pixels from datasets having $3200 \times 2048$ points, or possibly more. Therefore, implementation of $\Back$ represent the main 
bottleneck of the reconstruction process. If certain useful mathematical properties of $\Back$  are exploited, the computational
effort can be significantly reduced to $O( N^2 \log N )$ flops \cite{fftw,marone}.

Several techniques were developed aiming at a reduction to $O( N^2\log N )$ flops for computing $\Back$. One approach was established in \cite{andersson}, where the computation of $\Back g$ 
is performed after a change from cartesian to log-polar coordinates in the data. This approach leads to a convolution, which is computable through Fast Fourier Transform (\textsc{fft}) algorithms. Although elegant, the methods suffer from the ill-conditioning of the Log-polar 
transform at the ``fovea''. Nevertheless, it is possible to translate the fovea to different regions of the cartesian plane, in
order to enclose the reconstruction region. This leads to the concept of \emph{partial-backprojection} which can be easily implemented
in a parallel form. Other methods for fast computation of $\Back$ were presented in \cite{brandt,bresler2,bresler1, miqueles},
using  a divide and conquer strategy based on hierarchical decompositions of the full backprojection which are simpler than the full backprojection. Hierarchical decompositions can be created both in the image~\cite{bresler1} or in the data~\cite{brandt,bresler2}. Yet another approach is based on Non-uniform Fast Fourier Transform (\textsc{nfft}) algorithms (see \cite{bresler1,potts,karsten}) and the so-called Fourier Slice Theorem (\cite{kak_slaney, natterer}. In this paper, we propose another fast method for the computation of $\Back g$, also based on Fourier transforms. We claim that the backprojection of $g \in V$ can be easily done by filtering the lines of the $\tilde{g}$ one by one, where $\tilde{g}$ is the polar representation of $g$ in $S_+ = \R \times [0,2\pi]$.

This manuscript is organized as follows: Section \ref{sec:low} presents a discussion of low-complexity algorithms for the computation of $\Back$. Our 
low-complexity formula is presented in Section \ref{sec:our} and a discussion of the implementation is presented in Section \ref{sec:algo}.
Further comparison of all algorithms is presented in Section \ref{sec:all} and a discussion of the results is shown in Section \ref{sec:conc}.

\bigskip

\noindent{\bf Remark:} In this manuscript we use the the integral operator, sometimes with $\mathrm d x$ 
placed before the integrand, as it is more convenient to make explicit the variables being considered.
Whenever the integrand is short, we adopt the classic notation $\int f(x) \mathrm dx$.


\section{Class of Algorithms for the Backprojection}
\label{sec:low}

Let $g \in V$ be a given sinogram. A na\"ive implementation of the typical backprojection formula (\ref{adjoint}) has to be done using nested loops. Indeed, for each pixel $\bb x$ lying on a predefined meshgrid 
within the square $\|\bb x\|_\infty \leq 1$, the approximation of $b(\bb x) = \Back g(\bb x)$ is given by
\begin{equation} \label{eq:bp_app}
b(\bb x) \approx \Delta\theta \sum_{k=1}^{N_\theta} g( \bb x \cdot \bb \xi_{\theta_k}, \theta_k).
\end{equation}
It is easy to realize that the above approximation has a computational cost of $O(N_\theta)$ for each pixel $\bb x$, 
where $N_\theta$ is the total number of sampled angles. For a high resolution frame (see Figure \ref{fig:frame}.A), a linear interpolation for $\bb x \cdot \bb \xi_\theta$ on the grid of $-1 \leq t \leq 1$ is usually precise enough. Assuming that $b$ is represented by a square image of order $N \times N$, the total cost for computing the final backprojected image $b$ is $O(N^2 N_\theta)$. In practice, $N_\theta$ has almost the same magnitude of $N$, and thus we can state that the asymptotic cost to obtain $b$ is $O(N^3)$. Such an algorithm is impractical for high-resolution images.

There are at least three other types of backprojection algorithms which can dramatically reduce the computing time of the backprojected image $b$, for large datasets:
\begin{description}
\item[(i)] A fast slant-stack based approach \cite{fslant} was proposed by Averbuch {\it et al}. 
	Although this is an elegant and fast approach, it will not be covered in this manuscript;
\item[(ii)] \emph{Hierarchical decomposition} \cite{brandt,bresler1,bresler2}: Two different approaches that apply the divide-and conquer paradigm to the 
backprojection computation, splitting it recursively into smaller and simpler subproblems;

\item[(iii)] \emph{NFFT} \cite{nfft}: The Fourier Slice Theorem sets the Fourier Transform as bridge between the Radon 
	Transform $\mathcal R f( \theta, t )$ and the original image $f$. However, tomographic data does not provide 
	an evenly distributed sampling of the Fourier space, as required by traditional
	\textsc{fft} techniques (see \cite{cooley}). Use of this Fourier approach 
	was enabled by research on \textsc{nfft} algorithms (see \cite{bresler1,potts,karsten});
\item[(iv)] \emph{Anderson's formula} \cite{andersson}: Such a formula is obtained with an appropriate change of variables on the
	classical equation of the backprojection formula (\ref{adjoint}). The main idea is to convolve the sinogram in log-polar
	coordinates with an ideal kernel using \textsc{fft} algorithms.
\end{description}
In this paper, we focus mainly on the description of our algorithm and algorithm (iii).

\subsection{Log-polar backprojection}
\label{sec:and}

A fast method for obtaining the backprojection was derived by Andersson \cite{andersson}. His approach is based on a 
representation of the Backprojection/Radon transform as a convolution, by casting the computation in a log-polar 
coordinate system. In this section we propose a different proof for his formula.

Let $g = \Radon f$ be a given sinogram, i.e., the Radon transform of a compactly supported function $f$. 
Using the coordinate system notation of the previous section, where 
$[\cdot]_{\sf L}$ denotes the log-polar representation of some function, the main formula of the log-polar 
backprojection is written as
\begin{equation}\label{and_approach}
	[\Back g]_{\mathsf L}(\rho, \theta) = [g]_{\mathsf L} * [K]_{\mathsf L} \ (\rho, \theta),
\end{equation}
where $*$ stands for the two-dimensional convolution, and $K$ is the convolution Kernel
\begin{equation}\label{and_kernel}
	[K]_{\mathsf L}(\rho, \theta) = \delta(1 - e^{\rho}\cos\theta). 	
\end{equation}
Using above formula and the convolution theorem, we obtain
\begin{equation}
	[\Back g]_{\mathsf L} = \FT^{-1}\big( \FT [g]_{\mathsf L} \cdot \FT [K]_{\mathsf L} \big).
\end{equation}

Let us give a simple proof of the above equation, assuming that $f$ lies in a Schwartz space $\mathcal S(\R^2)$,
and $g \in \mathcal S(\R_+ \times [-\pi,\pi])$. 

\medskip

\noindent{\bf \it Proof:} We start with the integral representation of the 
backprojection operator, given in (\ref{eq:backCirc}) (See Appendix \ref{app:int}). Now, formula (\ref{and_approach}) is derived in four steps:
\begin{description}
\item[(a)] Changing the integral (\ref{eq:backCirc}) from cartesian coordinates $\bm{y}\in\R^2$ to Pr\"ufer coordinates, i.e., $\bm{y} \equiv \bm{y}_{\mu,\theta} =  p(\mu)\bm{\xi}_\phi$ we get $\mathrm d\bm{y} = |p'(\mu)p(\mu)| \mathrm d\mu \mathrm d\theta$ and
  \begin{equation} \label{eq:backPrufer}
    \Back g(\bm{x}_{\rho,\theta}) = \int_{S_+} g(\bm{y}_{\mu,\phi}) \delta \left( \kappa_{\bm{x}_{\rho,\theta}}(\bm{y}_{\mu,\phi}) \right) |p'(\mu)p(\mu)|\mathrm d \mu \mathrm d\phi
  \end{equation}

\item[(b)] The support of the Delta distribution in (\ref{eq:backPrufer}) is 
  \begin{equation} \label{eq:deltaPrufer}
    \kappa_{\bm{x}_{\rho,\theta}}(\bm{y}_{\mu,\phi}) = p(\mu)^2\left[ 1 - \frac{p(\rho)}{p(\mu)} \bm{\xi}_\phi \cdot \bm{\xi}_\theta  \right] = 
    p(\mu)^2\left[1- \frac{p(\rho)}{p(\mu)} \cos(\phi-\theta)\right]
  \end{equation}

\item[(c)] Let $[\cdot]_{\mathsf{Pr}}$ be the representation in Pr\"ufer coordinates. From (\ref{eq:deltaPrufer}) and
  (\ref{eq:backPrufer}) we arrive at
  \begin{multline}
    [\Back g]_{\mathsf{Pr}}(\rho,\theta)\\ {}= \ds\int_{S_+} [g]_{\mathsf G}(\mu,\phi) \delta \left(         p(\mu)^2\left[1- \frac{p(\rho)}{p(\mu)} \cos(\phi-\theta)\right] \right) |p'(\mu)p(\mu)| \mathrm     d \mu \mathrm d\phi \\
    {} = \ds \int_{S_+} [g]_{\mathsf G}(\mu,\phi) \delta \left( 1- \frac{p(\rho)}{p(\mu)}             \cos(\phi-\theta)\right) \frac{|p'(\mu)p(\mu)|}{p(\mu)^2} \mathrm d \mu \mathrm d\phi                \label{eq:quasiConv}
    \end{multline}
  where $S_+ = \R_+ \times [-\pi,\pi]$

\item[(d)] A convolution is obtained in (\ref{eq:quasiConv}) only if $p$ is such that $p(\rho)=p(\mu)p(\rho-\mu)$, which in turn implies that $p$ is an exponential function. Hence, Pr\"ufer coordinates reduce to log-polar coordinates, which we denote by $[\cdot]_{\mathsf L}$. Finally, we obtain
    \begin{equation} \label{eq:andersson}
      [\Back g]_{\mathsf L}(\rho,\theta) = \int_{S_+} [g]_{\mathsf L}(\mu,\phi)\delta\left(1-e^{\rho-\mu}\cos(\phi-\theta)\right) \mathrm d\mu
      \mathrm d\phi
    \end{equation}
which is the final convolution formula. $\qed$
\end{description}


\section{Back-projection Slice Theorem}
\label{sec:our}

Although Anderson's approach is asymptotically fast, it has a few drawbacks. Firstly, the gain of speed using Fourier transforms to compute the convolution 
is reduced with forward/backward log-polar transformations. Also, these interpolations can produce errors, especially near the origin, due to a strong non-uniformity of the log-polar mesh in that region. To avoid these factors, another approach for the calculation of the backprojection operator can be used. This approach is based on the following theorem: 

\begin{theo*}[\bf Backprojection Slice Theorem (BST)]
	 Let $g = g(t,\theta) \in V$ a given sinogram and\/ 
	 $\widehat{\cdot}$ denotes the Fourier transform operation. It follows that the backprojection $\Back$ satisfies 
	\begin{equation} \label{eq:bst}
	\widehat{\Back g}(\sigma\bb \xi_\theta) = \frac{\hat{g}(\sigma,\theta)}{\sigma}
	\end{equation}
	with $\sigma>0 \in \R$ and $\theta\in[0,2\pi]$.
\end{theo*}

\bigskip

\begin{proof}
Using the sifting property of the $\delta$-distribution, the backprojection (\ref{B}) can be presented 
in the following form
\begin{equation}\nonumber
\Back g(\bb x) = \int_0^{\pi} g( \bb x \cdot \bb \xi_\theta,\theta) \mathrm d \theta 
= \int_0^{\pi} \int_\R g(t,\theta) \delta (t- \bb x\cdot\bb \xi_\theta) \mathrm d t \mathrm d \theta 
\end{equation}
Considering the two-dimensional Fourier transform of $\Back g$, i.e., $\FT\colon \Back g \mapsto \widehat{\Back g}$ and using representation (\ref{eq:backCirc}) (see Appendix),
\begin{eqnarray*}
\widehat{\Back g}(\bb\omega) = \int_{\R^2} \Back g(\bb x) e^{-i \bb \omega\cdot \bb x} \mathrm d\bb x
= \int_{\R^2}\!\!\! \mathrm d\bb x \int_{\R^2}\!\!\! \mathrm d\bb y \ [g]_{\mathsf c}(\bb y) \delta \left(\bb y \cdot(\bb y - \bb x)\right)
e^{-i\bb \omega\cdot \bb x}  \\
= \int_{\R^2}\!\!\! \mathrm d \bb y \  [g]_{\mathsf c}(\bb y) \int_{\R^2}\!\!\! \mathrm d\bb x \ \delta \left(\bb y\cdot(\bb y-\bb x)\right) e^{-i\bb \omega\cdot \bb x}
\equiv  \int_{\R^2}\!\!\! \mathrm d\bb y \  [g]_{\mathsf c}(\bb y) \mathcal{T}(\bb y,\bb \omega)
\end{eqnarray*}
where $\bb y, \bb \omega \in \R^2$ and 
\begin{equation} \label{GGG}
	\mathcal{T}(\bb y,\bb \omega) = \int_{\R^2}\!\!\!\mathrm d\bb x \ \delta \left( h_{\bb y}(\bb x)\right) e^{-i\bb \omega\cdot\bb x}, \ \ \ \ h_{\bb y}(\bb x) = \bb y \cdot(\bb y-\bb x)
\end{equation}

Since the distribution (\ref{GGG}) is supported in the set $h^{-1}_{\bb y}(0) = \{\bb x\in\R^2 \colon h_{\bb y}(\bb x) = 0\}$, 
it follows from (\ref{eq:pathint}) (See Appendix \ref{app:int}) and $\nabla h_{\bb y} = -\bb y$ that
\begin{equation} \label{eq:T}
\mathcal{T}(\bb y,\bb \omega) = \frac{1}{\|\bb{y}\|} \int_{h^{-1}_{\bb{y}}(0)} e^{-i \bb{\omega}\cdot \bb{x}} \mathrm d s(\bb x)
= \int_{h^{-1}_{\bb{y}}(0)} e^{-i \bb{\omega}\cdot \bb{x}(q)} \mathrm dq
\end{equation}
The set $h^{-1}_{\bb{y}}(0)$ determines a straight line passing through $\bb{y}$ and with normal vector $\bb{y}$. Thus, 
$h^{-1}_{\bb{y}}(0) = \bb{y} + \mbox{span}\{S\bb{y}\}$, being $S\bb{y} \perp \bb{y}$ and $S$ a $\frac{\pi}{2}$-rotation
matrix. Therefore, $\bb{x}(q) \in h^{-1}_{\bb{x}}(0)$ is on the form $\bb{x}(q) = \bb{y} + q S\bb{y}$ and the integral
in (\ref{eq:T}) can be written as:
\begin{equation} 
\mathcal{T}(\bb{y},\bb{\omega}) = \int_{\R} e^{-i \bb{\omega}\cdot [\bb{y}+q S\bb{y}]} \mathrm d q 
=  e^{-i \bb{\omega}\cdot \bb{y}} \int_{\R} e^{-i q \bb{\omega}\cdot (S\bb{y})} \mathrm d q = e^{-i \bb{\omega}\cdot \bb{y}} \delta\left(\bb{\omega}\cdot S\bb{y} \right)
\end{equation}
Hence, the Fourier transform of $\Back g$ becomes
\begin{equation}  \label{eq:quasi2}
\widehat{\Back g}(\bb{\omega}) = \int_{\R^2} [g]_{\mathsf c}(\bb{y}) \delta\left(\bb{\omega}\cdot S\bb{y} \right) e^{-i \bb{\omega}\cdot \bb{y}} \ \mathrm d\bb{y} 
\end{equation}
For $\bb{\omega}$ fixed, $\{\bb{y}\in\R^2 \colon \bb{\omega}\cdot (S\bb{y}) = 0\} = \mbox{span}\{\bb{\omega}\}$, with $S$ aπ $\frac{\pi}{2}$-rotation matrix. Indeed, since $S\bb{y} \perp \bb{w}$ and $S\bb{y}\perp \bb{y}$, it follows $\bb{\omega}\parallel \bb{y}$. Once again, using the representation
(\ref{eq:pathint}) (see Appendix \ref{app:int}) for  (\ref{eq:quasi2}) we arrive at
\begin{equation}  \label{eq:quasi3}
\widehat{\Back g}(\bb{\omega}) = \int_{\R} \frac{ [g]_{\mathsf c}(q\bb{\omega})}{\|S\bb{\omega}\|} e^{-i  \bb{\omega}\cdot(q\bb{\omega})} \ \mathrm ds(\bb \omega) 
\end{equation}
Since $\|S\bb{\omega}\| = \|\bb{\omega}\|$ and $\mathrm d s(\bb\omega) = \|\bb \omega\| dq$, we finally obtain
\begin{equation} \label{eq:FourierBack1}
\widehat{\Back g}(\bb{\omega}) = \int_{\R} [g]_{\mathsf c}(q\bb{\omega}) e^{-i q \|\bb{\omega}\|^2} \ \mathrm dq
\end{equation}
From the above equation, the backprojection is a polar convolution. Indeed, switching the frequency
domain to polar coordinates, i.e., $\bb{\omega} = \sigma \bb{\xi}_{\theta}$ (with $\sigma\in\R_+$ and $\theta \in [0,2\pi]$) we get 
\begin{equation} \label{eq:J}
\widehat{\Back g}(\sigma\bb{\xi}_\theta) = \int_{\R} [g]_{\mathsf c}(q\sigma\bb{\xi}_\theta) e^{-i q \|\sigma\bb{\xi}_\theta\|^2} \mathrm dq
= \int_{\R} \frac{[g]_{\mathsf c}(u\bb{\xi}_\theta)}{\sigma} e^{-i u \sigma} \mathrm du.
\end{equation}
Now, letting $[\cdot]_{\mathsf s}$ be the representation in semi-polar coordinates, it is true that
$[g]_{\mathsf c}(u\bb{\xi}_\theta) = g(u,\theta)$ is the input sinogram $g(u,\theta)$. 
From (\ref{eq:FourierBack1}) and (\ref{eq:J}), using polar coordinates
\begin{equation}\label{eq:sliceBack}
[\widehat{\mathcal \Back g}]_{\sf p}(\sigma,\theta) = \widehat{\mathcal \Back g}(\sigma\bb{\xi}_\theta) = \frac{1}{\sigma} \int_{\R} g(u,\theta) e^{-i u \sigma} \mathrm du
\end{equation}
Identity (\ref{eq:sliceBack}) is our backprojection-slice Theorem (\ref{eq:bst}) for computing the operator $\Back$.
\end{proof}

Indeed, at each radial line $\theta$ in the frequency domain, the two-dimensional Fourier transform of $\Back$
equals the one-dimensional radial Fourier transform of the projection $g(t,\theta)$ multiplied by the kernel $1/\sigma$ for $\sigma>0$.

\medskip

\noindent {\bf Remark 1:} The mathematical proof outlined above provides a \emph{direct} formula for
the computation of a backprojected image, i.e., given a sinogram $g$, the explicit steps to compute 
the backprojection in the frequency polar coordinates results in formula (\ref{eq:bst}). In practice, 
there are several iterative methods that depends explicitly on the computation of the backprojection 
of \emph{any} sinogram. In the other hand, analytical formulas usually handle with the backprojection 
of a filtered sinogram, from where standard formulas like the \emph{filtered backprojection} or
the \emph{filter of the backprojection} are established. To validate our backprojection result
we remark the following items:

\medskip
\noindent{\bf (i)} It is a well known fact \cite{deans,helgason, kak_slaney} that, for a given feature function $f \in U$,
the following property holds
\begin{equation} \label{eq:BRf}
\Back \Radon f (\bb x) = (f * h)(\bb x), \ \ \ \ h(\bb x) = \frac{1}{\|\bb x\|}
\end{equation}
which, in the frequency domain, is written as (cartesian and polar representation, respectively)
\begin{equation} \label{eq:fob}
\widehat{\Back \Radon f}(\bb w) = \hat{f}(\bb w) \frac{1}{\|\bb w\|} \ \ \ \Leftrightarrow \ \ \ 
\widehat{\Back \Radon f}(\bb \sigma \bb \xi_\theta) = \hat{f}(\bb \sigma \bb \xi_\theta) \frac{1}{\sigma}  
\end{equation}
due to the fact that $\FT \colon \frac{1}{\|\bb x\|} \mapsto \frac{1}{\|\bb w\|}$. Now, replacing the backprojection
slice theorem (\ref{eq:bst}) into (\ref{eq:fob}), we obtain
\begin{equation}
\widehat{\Radon f}(\sigma \bb \xi_\theta)\frac{1}{\sigma} = \hat{f}(\bb \sigma \bb \xi_\theta) \frac{1}{\sigma} \ \ \ \Rightarrow \ \ \ 
\widehat{\Radon f}(\sigma \bb \xi_\theta) = \hat{f}(\bb \sigma \bb \xi_\theta)
\end{equation}
which is the celebrated Fourier Slice-Theorem \cite{nattw}. 

\medskip
\noindent{\bf (ii)} From the classical inversion of the Radon transform, i.e., the filtered-backprojection
algorithm, it is true that
\begin{equation} \label{eq:fbp}
\Back F g (\bb x) = f(\bb x), \ \ \  g = \Radon f
\end{equation}
where $F$ is a low-pass filtering operator, that is $\widehat{Fg}(\nu,\theta) = \hat{g}(\nu,\theta) |\nu|$, for $\nu \in \R$. In the polar frequency domain, (\ref{eq:fbp}) reads $\widehat{\Back Fg}(\sigma \bb\xi_\theta) = \hat{f}(\sigma \bb\xi_\theta)$. From the backprojection slice Theorem (\ref{eq:bst}), such equation becomes
\begin{equation}
\frac{1}{\sigma}\widehat{F g}(\sigma \bb\xi_\theta) = \hat{f}( \sigma \bb \xi_\theta) \ \ \ 
\Rightarrow \ \ \ \frac{1}{\sigma}\hat{g}(\sigma \bb\xi_\theta) \sigma = \hat{f}(\sigma \bb \xi_\theta), \ \ \ \ \sigma \in \R_+
\end{equation}
Once again, the above equation yields the Fourier Slice-Theorem.

\bigskip

\noindent {\bf Remark 2:} The \textsc{dc}-component of the Backprojection of some function $g$ lying in the sinogram space is
defined by
\begin{eqnarray}
\widehat{\Back g}(\bb 0) &=& \int_{\R^2} \Back g(\bb x) \mathrm d \bb x \\
&=& \int_{\R^2} \mathrm d \bb x \int_0^\pi \mathrm d\theta \ g(\bb x\cdot \bb \xi_\theta, \theta) \\
&=& \int_0^\pi \mathrm d \theta \int_{\R} \mathrm d t \int_\R \mathrm d s \ g(t, \theta) \equiv M 
\end{eqnarray}
where we have used $\mathrm d\bb x = \mathrm d t\mathrm d s$ to make explicit the change of
variables from $\bb x$ to $(t,s)$, being $s$ the variable along the direction $\bb \xi_\theta^\perp$.
The \textsc{dc} of an arbitrary $g \in V$ provide\footnote{Even if $g$ is the sinogram of a compactly supported function feature function on the unit disk $\|x\|_2 \leq 1$, we have $\hat{g}(0, \theta) = \mbox{constant}$, although with
$M = \infty$.} $M = \infty$. In this sense, $\widehat{\Back g}$ behaves like a tempered
distribution since $\Back g$ lies in a Schwartz space, where the Fourier transform is an automorphism. 
Also, it is easy to note that
\begin{equation} \label{eq:ratio}
\frac{\hat{g}(\sigma,\theta)}{\sigma} = i \hat{h}(\sigma,\theta), \ \ h(t) = \int_{\infty}^t g(t,\theta) \mathrm d t 
\end{equation}
i.e., $h$ is a primitive of $g$. Hence, using (\ref{eq:ratio}) as $\sigma \to 0$, the limit of the ratio 
$\hat{g}(\sigma,\theta)/\sigma$ diverge in $\sigma = 0$. Finally, \textsc{bst} formula can be easily applied for
some $g \in V$ with a nonzero \textsc{dc}-component. In fact, setting $p(t,\theta) = g(t,\theta) - \hat{g}(0,\theta)$, 
it is true that $\hat{p}(0,\theta)=0$ and the backprojection of $g$ follows with $\Back g(\bm x) = \Back p(\bm x) + \hat{g}(0,\theta)$.


\section{Analytical Examples}
\label{sec:anal} 

In this section, we provide two examples where the analytical computation of the 
backprojection operator is possible. This is important to validate further numerical simulations.
The first example given is for a point source function, both in log-polar and 
cartesian coordinates. The second one, for a symmetrical circular function. In 
what follows, we consider that $g = g(t,\theta)$ is the Radon transform of 
$f = f(\bm x)$, while $b = b(\bm x)$ is the final backprojected image.

\bigskip

\noindent{\bf Example I:} Following Andersson's formula \label{eq:andersson}, the backprojection of any sinogram $g$ is written
as a convolution in Log-polar coordinates
\begin{equation} \label{eq:and}
b(e^\rho \bm \xi_\theta) = \Back g( e^{\rho} \bm \xi_\theta) = \int \mathrm d u \int \mathrm d\beta \ \ g(e^u, \beta) \delta (1 - e^{\rho-u} \cos (\theta-\beta) )
\end{equation}
It is a well known fact that the Radon transform of a single point source, located 
at $\bm x= \bm a$ is
\begin{equation} \label{eq:data}
f(\bm x) = \delta(\bm x-\bm a) \ \ \ \implica \ \ \ g(t,\theta) = \delta(t-\bm a \cdot \bm \xi_\theta)
\end{equation}
Taking $\bm a = e^A \bm \xi_\phi$ as the log-polar representation of the
source point $\bm a$, we use (\ref{eq:data}) and (\ref{eq:and}) to obtain $b$ as
\begin{eqnarray}
b(e^\rho \bm \xi_\theta) &=& \int \mathrm d \beta \int \mathrm d u \ \delta( \bm a\cdot\bm \xi_\beta - e^u ) \delta( 1 - e^{\rho-u} \cos (\theta-\beta)) \\
&=& \frac{1}{\sqrt{(\cos(\theta-\phi) e^\rho)^2 + (e^A - \sin(\theta-\phi) e^\rho)^2}} \label{eq:bplogp}
\end{eqnarray}
where $A = \ln \|a\|$. The details of the log-polar representation (\ref{eq:bplogp}) 
are presented in the Appendix \ref{sec:psf}. To obtain a cartesian representation 
we use $\bm x = e^\rho \bm \xi_\theta$ and
\begin{equation}
\cos(\theta-\phi) = \bm \xi_\theta\cdot \bm \xi_\phi, \ \ \ \ \
\sin(\theta-\phi) = \bm \xi_\theta \cdot \bm \xi_\phi^\perp
\end{equation}
Now, (\ref{eq:bplogp}) becomes
\begin{eqnarray}
b(\bm x) &=& \frac{1}{\sqrt{[(\bm \xi_\theta\cdot \bm \xi_\phi) e^\rho]^2 + [e^A - (\bm \xi_\theta \cdot \bm \xi_\phi^\perp) e^\rho]^2}} \\
&=& \frac{1}{\sqrt{[(e^\rho \bm \xi_\theta)\cdot (e^A\bm \xi_\phi) e^{-A}]^2 + [e^A - (e^\rho \bm \xi_\theta) \cdot (e^A \bm \xi_\phi^\perp) e^{-A}]^2}} \\ 
&=& \frac{1}{\sqrt{[(\bm x \cdot \bm a) e^{-A}]^2 + [e^A - (\bm x \cdot \hat{\bm a}) e^{-A}]^2}}
\end{eqnarray}
where $\hat{\bm a} = e^A \bm \xi_\phi^\perp$ is a counterclockwise rotation by $\frac{\pi}{2}$ of the point source $\bm a$. Finally, since $e^{A} = \|\bm a\|$ we obtain  
\begin{equation}
b(\bm x) = \frac{\|\bm a\|}{\sqrt{[\bm x \cdot \bm a ]^2 + [ \|\bm a\|^2 - \bm x \cdot \hat{\bm a}]^2}} \label{eq:bpcart}
\end{equation}
which is the cartesian representation of the backprojection of the sinogram $g$
given in (\ref{eq:data}).

\bigskip

\noindent{ \bf Example II:} Considering $f(\bm x) = \mbox{circ}(\|\bm x\|)$, i.e.,
\begin{equation}
f(\bm x) = \left\{\begin{array}{ll}
1, & \|\bm x\| \leq 1 \\
0, & \|\bm x\| > 1 \\
\frac{1}{2}, & \|\bm x\| = 1
\end{array} \right.
\end{equation}
it is known \cite{bracewell} that 
\begin{equation} \label{eq:Fcirc}
\hat{f}(\bm \omega) = \frac{J_1(\|\bm \omega\|)}{\|\bm \omega\|}
\end{equation}
with $J_1$ the order 1 Bessel function of the first kind. From
the Fourier-Slice-Theorem $\hat{f}(\sigma \bm \xi_\theta) = \hat{g}(\sigma,\theta)$ 
and (\ref{eq:Fcirc}) it is easy to obtain $g$ satisfying
\[
\hat{g}(\sigma,\theta) = \frac{J_1 (\|\sigma \bm \xi_\theta\|)}{\|\sigma \bm \xi_\theta\|}  = \frac{J_1(\sigma)}{\sigma}
\]
Now, using the \textsc{bst} formula (\ref{eq:bst}), the Fourier representation for $b$ 
becomes
\begin{equation} \label{eq:bp_circ}
\hat{b}(\sigma \bm \xi_\theta) = \frac{\hat{g}(\sigma,\theta)}{ \sigma} 
= \frac{J_1(\sigma)}{\sigma^2} \  \iff  \ 
\hat{b}(\bm \omega) = \frac{J_1(\bm \omega)}{\|\bm \omega\|^2}
\end{equation}
The above equation provides a testing algorithm (in the Fourier domain)
for our numerical strategies. In fact, either \textsc{bst} or Andersson's
algorithm presents difficulties as $\bm \omega \to \bm 0$, as discussed in 
the next section. 


\section{Implementation issues}
\label{sec:algo}

All the notation needed for the implementation of the \textsc{bst} formula is presented in Table \ref{tab:symbols}.
In this section we assume to have the sinogram presented at the nodes of the uniform {\it sinogram grid} 
$G_{\sf s}$. Quantification of the sinogram function $g(\theta,t)$ over $G_{\sf s}$ will be denoted by
$g_{ij}$.

\begin{table}[h]
\footnotesize
\begin{tabular}{l|lll}
\hline
\specialcell{Mesh\\sequence}  & \specialcell{Description\\of coordinates} & \specialcell{Mesh\\Size} & \specialcell{Step\\size} \\
\hline
$\{t_k\} = \{-1.0, \ldots, 1.0\}$ & sinogram & $|t_k| = N_t$ & $\Delta t = 2/N_t$ \\
\hline 
$\{s_k\} = \{0, \ldots, 1.0\}$ & polar & $|s_k| = N_s = N_t/2$ & $\Delta s = 2/N_t$ \\
\hline 
$\{\rho_k\} = \{\rho_0, \ldots, \ln 1\}$ & log-polar & $|\rho_k| = N_\rho$ & $\Delta \rho = (-\rho_0 + \ln 1)/N_\rho$ \\
\hline
$\{\theta_k\} = \{0, \ldots, \pi\}$ & angles within $[0,\pi]$ & $|\theta_k|=N_\theta$ & $\Delta \theta = \pi/N_\theta$  \\
\hline
$\{\phi_k\} = \{0, \ldots, 2\pi\}$ & angles within $[0,2\pi]$ & $|\phi_k|=2N_\theta$ & $\Delta \phi = \pi/N_\theta$  \\
\hline	
$G_{\sf s} = \{\theta_k\}\times \{t_i\}$ & 2D sinogram mesh & $|G_{\sf s}| = (N_t, N_\theta)$ & - \\
\hline
$G_{\sf P} = \{\phi_k\}\times \{s_i\}$ & 2D polar mesh & $|G_{\sf P}| = (N_t/2, 2N_\theta)$ & - \\
\hline
$G_{\sf L} = \{\phi_k\}\times \{\rho_i\}$ & 2D log-polar mesh & $|G_{\sf L}| = (N_\rho, N_\theta)$ & - \\
\hline
\end{tabular}
\caption{Glossary of symbols used for implementation details.}
\label{tab:symbols}
\end{table}

\subsection{Algorithm for Log-polar backprojection}

In the implementation of the simplest algorithm, based on the Andersson's approach, we don't take 
into account the irregularity in the origin. The algorithm can be summarized as follows: 

\begin{description}

\item[Step 1.] {\it Interpolate the sinogram to log-polar coordinates, $[g]_{\mathsf s} \to [g]_{\mathsf L}$}: 
The domain of the experimentally obtained sinogram is $S_+ = [0, \pi] \times [-1, 1]$. 
To translate it to Log-Polar coordinates it is more convenient to first translate sinogram coordinates to 
standard semi-polar coordinates (see Fig. \ref{fig:phantom_sino}). In this case, 
now we have the sinogram in semi-polar coordinates, sampled in the nodes of the {\it polar grid} $G_{\sf P}$
(see Table \ref{tab:symbols}). The change from semi-polar to log-polar coordinates can be easily done using the linear interpolation along the ray 
for every $\theta=Const$. Let us denote the number of points along every ray as $N_{\rho}$. Since $\ln(\epsilon) \to -\infty$ 
as $\epsilon \to 0$ (we will discuss about this disadvantage later), we have to select a number $\rho_0 < 0$ to define the lowest 
point of the mesh, closest to the origin. Thus, we have to interpolate sinogram from the mesh $G_{\sf s}$ to the {\it log-polar grid}
$G_{\sf L}$ (see Table \ref{tab:symbols}). In fact, for all change of coordinates computed in this work, we use simple linear interpolation. The problem of selection of the 
first node $\rho_0$ of log-polar mesh considered below, in Sec. \ref{sec:origin}.


One of the specific features of Log-Polar mesh is its non-uniformity. To obtain clear interpolation without losing information, we need to make the biggest step of log-polar mesh equal to radial step $\Delta s$ of the original sinogram. It
can be easily done by finding the mesh size $N_{\rho}$ from the equation:
\begin{equation}
	\exp(0) - \exp(-\Delta \rho) = \Delta s
\end{equation}
The above equation comes from the definition of the mesh $G_{\sf L}$ since $\{\rho_k,\rho_{k-1}\}$ is in fact 
$\{0,- \Delta \rho\}$ with $k = N_\rho - 1$. Since $\Delta \rho = \rho_0/N_{\rho}$ we finally obtain the number of points
in the log-polar system 
\begin{equation}\label{N_rho}
	N_{\rho} = - \frac{\rho_0}{\ln(1-\Delta s)}.
\end{equation}
In the next Section we consider a different choice for the parameter $\rho_0$.

\item[Step 2.] {\it Calculate the kernel $K$ by formula (\ref{and_kernel})}: The kernel 
represented with the formula (\ref{and_kernel}): $[K]_{\sf L}(\rho, \theta) = \delta(1 - e^{\rho}\cos\theta)$ and can be approximated on the mesh 
$G_{\sf L}$ using the condition: 
\begin{equation}\label{calc_and_kernel}
	k_{ji} = \begin{cases} 
			1/\Delta\rho, & \quad \text{if} \quad | e^{\rho_i}\cos\theta_j - 1 | \leq \Delta \rho, \\
			0, & \quad \text{otherwise}. 
		\end{cases}
\end{equation}
Here $\Delta\rho = -\frac{\rho_0}{N_{\rho}}$ is the step on the mesh for variable $\rho$. The kernel in log-polar 
coordinates and the absolute value of its Fourier transform are shown on Figure \ref{fig:and_kernel}. The 
approximation (\ref{calc_and_kernel}) could be numerically improved using appropriate strategies for 
the evaluation of a Delta distribution concentrated on the zero level set of a function, see \cite{towers}.

\begin{figure}[t!]
    \centering
    \begin{tabular}{cc} 
    (a) & (b) \\
    \begin{tikzpicture}
        \begin{axis}[width=6cm, compat=newest,
        ticks=none,enlargelimits=false, axis on top, scale only axis, 
        axis equal image, xlabel={$\theta \in [-\pi,\pi]$}, ylabel={$\rho \in [\rho_0,0]$}]
        \addplot graphics 
        [xmin=0,xmax=256,ymin=0,ymax=256]{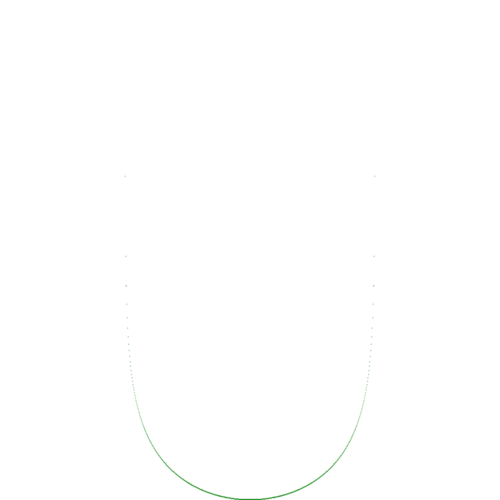};
        \end{axis}
    \end{tikzpicture} &
    \begin{tikzpicture} 
        \begin{axis}[width=6cm, compat=newest,
        ticks=none,enlargelimits=false, axis on top, scale only axis, 
        axis equal image,xlabel={$K(\rho,\theta)$}]
        \addplot graphics 
        [xmin=0,xmax=256,ymin=0,ymax=256]{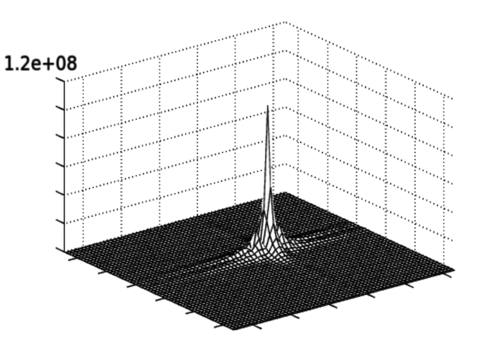};
        \end{axis}
    \end{tikzpicture}
    \end{tabular}
    \caption{\small Kernel in log-polar representation (a) and its Fourier image (b)}
    \label{fig:and_kernel}
\end{figure}

\item[Step 3.]{\it Calculate the convolution $[K]_{\mathsf L}\star[g]_{\mathsf L}$}: Using the uniform grids, 
the convolution is calculated using the Fast Fourier Transform pair through \textsc{fftw3} software library.

\item[Step 4.]{\it Interpolate the result from previous step back to Cartesian coordinate system}. We are 
using bilinear interpolation, which already has been described in Step 1.

\end{description}

\subsection{Problem near the Origin}\label{sec:origin}

Due to the log-polar representation, $s \to 0$ causes $\rho \to -\infty$. Of course, in real 
calculations it is not possible to get a proper interpolation to this grid. This practical problem can be 
solved in two ways. The first approach is clearly mathematical, and was proposed by Andersson in his work 
\cite{andersson}, called \emph{partial backprojection method}. This method is based on 
moving the origin outside the region of interest, which allows us to make a clear 
interpolation of all points of the sinogram with non-zero values. The second approach is to select 
a proper $\rho_0$, which adapts properly to the resulting cartesian grid. As we will show, 
this way also gives good results. 

\paragraph{\textbf{Adaptive selection of $\rho_0$.}}

For clear interpolation we have to use $\rho_0$ as a very big negative number, from which we start the approximation
to the log-polar mesh. But, in fact, this number is connected to the mesh which is chosen for cartesian
representation of the result. Assume the cartesian mesh with the direction steps $(\Delta x, \Delta y)$.
In this case, to avoid the loss of information near the origin, we have to set $\rho_0 < \ln(\min(\Delta x, \Delta y))$. 
The result of using the Anderson's with the rightly chosen $\rho$ is presented Fig.\ref{fig:logpolar_origin}.(a) in the ``Numerical results'' section.
In fact, it is important to note that the performance of Log-polar backprojection depends on the desired 
resolution of the resulting image. The oversampling of log-polar mesh grows up fast as the number of pixel 
increases in the cartesian grid, as it shown below in this section. 

Now let $g(t,\theta)$ be some sinogram and $[g]_{\sf L}(\rho, \theta)$ his
log-polar representation for $\rho\in(-\infty, 0)$. 
Consider the approximation of $[g]_{\sf L}$ with the following compactly supported function: 
\begin{equation}\label{gl_approx}
	g_{\sf L}^{-}(\rho, \theta) = 
	\begin{cases}
		[g]_{\sf L}(\rho, \theta), & \rho \in [\rho_0, 0], \\
		0, & \rho < \rho_0  
	\end{cases}
\end{equation}
where $\rho_0 \ll 0$ is a given fixed parameter. We want to measure the 
norm of discrepancy between $\Back g$ and $\Back g^-$, i.e., 
\begin{equation}
\begin{array}{lll}
	\| [\Back g]_{\sf L} - [\Back g^{-}]_{\sf L} \|_{L_2}^2 = \int\limits_{0}^{2\pi}\int\limits_{-\infty}^{0} 
		\big(  [\Back g]_{\sf L}(\rho, \theta) - [\Back g^{-}]_{\sf L}(\rho, \theta) \big)^2 e^{2\rho}\mathrm d\rho \mathrm d\theta \\
		\ \ \ = \displaystyle \int\limits_{0}^{2\pi}\int\limits_{\rho_0}^{0} + \int\limits_{0}^{2\pi}\int\limits_{-\infty}^{\rho_0} \ 
		\big(  [\Back g]_{\sf L}(\rho, \theta) - [\Back g^{-}]_{\sf L}(\rho, \theta) \big)^2 e^{2\rho}\mathrm d\rho \mathrm d\theta
\end{array}
\end{equation}
From (\ref{gl_approx}) it follows that
\begin{equation} \label{eq:error}
	\| [\Back g]_{\sf L} - [\Back g^{-}]_{\sf L} \|_{L_2}^2 = 
	\int\limits_{0}^{2\pi}\int\limits_{-\infty}^{\rho_0} 
	[\Back g]_{\sf L}(\rho, \theta)^2 e^{2\rho}\mathrm d\rho \mathrm d\theta \le 2\pi {\sf c}^2 e^{2\rho_0},
\end{equation}
where $\sf c = \displaystyle\max_{\rho \le \rho_0} [\Back g]_{\sf L}$ is a constant, which, in practice, refers to the value of Backprojection in the origin. 
This value can be easily estimated. Equation (\ref{eq:error}) give us a bound 
for the error, when we remove the origin in the computation of the backprojected image.

To obtain a good reconstruction, it is easy to obtain the number $N_\rho$ from using the following formula 
\begin{equation}\label{N_rho}
	N_{\rho} \approx \frac{\ln(\min(1/N_x, 1/N_y)}{\ln(1 - \Delta s)}.
\end{equation}
For example, assuming that $N_s = 1024$ and $N_x = N_y = 1024$, then $N_{\rho} = 3546$. Therefore, an oversampling 
of the input data is usually needed (about $\approx 4$ times, which can be easily estimated from formula \ref{N_rho}). This fact, of course, 
decreases the calculation speed of convolution. In our implementation, we have obtained higher speed for 
partial backprojections (described below) due to few interpolation steps. The processing speed of these backprojection formulas
are considered in Section \ref{sec:all}.

\paragraph{\textbf{Partial Backprojections.}}

\begin{figure}[ht!]
	\begin{tabular}{cc}
		(a) & (b) \\
		\begin{tikzpicture} 
	        \begin{axis}[width=6cm, compat=newest,
	        ticks=none,enlargelimits=false, axis on top, scale only axis, 
	        xlabel={$\bm x_1$}, ylabel={$\bm x_2$}]
	        \addplot graphics 
	        [xmin=0,xmax=256,ymin=0,ymax=256]{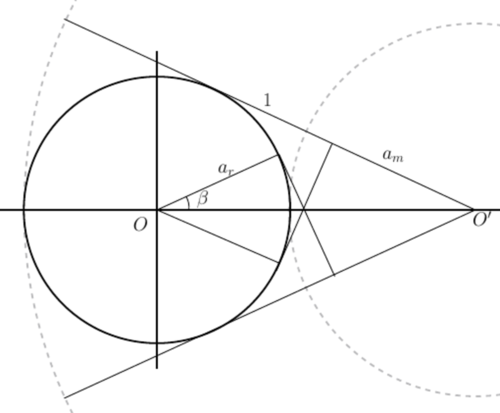};
	        \end{axis}
	 	\end{tikzpicture}
		 & 
		\begin{tikzpicture} 
	        \begin{axis}[width=6cm, compat=newest,
	        ticks=none,enlargelimits=false, axis on top, scale only axis, 
	        axis equal image, xlabel={$\theta\in[-\pi,\pi]$}, ylabel={$\rho \in [\rho_0,0]$}]
	        \addplot graphics 
	        [xmin=0,xmax=256,ymin=0,ymax=256]{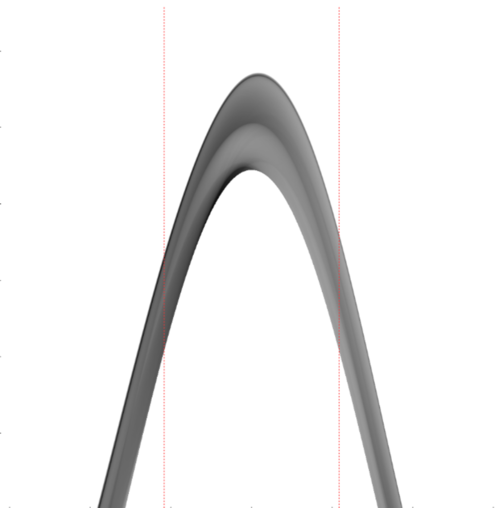};
	        \end{axis}
	 	\end{tikzpicture}
	\end{tabular}
	\caption{\small (a) The scheme of the reconstruction with sectoral method; 
		(b) transformed sinogram and sector outside of zero (inside red lines)}
	\label{fig:sectoral_scheme}
\end{figure}

The Partial Backprojections method is based on the \emph{shifting property} of the Radon
transform, defined as
\begin{equation}\label{sh_prop}
	u(\bm x) = f(\bm x - \Delta \bm x)\ \ \implica \ \ \Radon u(t,\theta) = g(t - \bm \xi_\theta \cdot \Delta \bm x, \theta), \ \ \ g= \Radon f
\end{equation}
Using this formula we can transform the original sinogram, presented in semi-polar coordinates 
and take the part, which is located far away from zero. We choose some angle $\beta$ and consider 
the sector of the original sinogram $\theta \in [\theta_0, \theta_1]$, where $\theta_1 = \theta_0 + 2\beta$. 
Rescaling the original sinogram to the size $a_r$, as it shown on Fig.\ref{fig:sectoral_scheme}.(a) and 
rotating it in the way that sector under investigation will be located in $\theta \in [-\beta, \beta]$.
Now, it is easy to obtain the values of the distances between old and new origins $(1-a_r)$ and between 
the sector under investigation and new origin $1-2a_r$, which defines the minimal $s$ and $\rho_0$ for 
the interpolation from semi-polar to log-polar coordinates. This method is described in details
in \cite{andersson}. 

Using partial backprojection is convenient because it is possible to highlight any sector of interest of the sinogram without any information loss. Also, it is not necessary to process all the sinogram at once - 
we can process only the parts that we are interested in. The first disadvantage of this method is that the Fourier image of the kernel is singular at the origin (see Fig.\ref{fig:and_kernel}.(b)), which causes artifacts on the result 
of sector backprojection. We note that, mathematically, this problem also exists in the previous method (adaptive choosing 
the $\rho_0$), but, since we exclude the origin from the calculations, we just "skip" this singularity. 
The second disadvantage of the partial method is that further mathematical operations are needed, e.g., translation 
of the origin in different coordinate systems (log-polar and Cartesian) and an extension of angles of the sectors under reconstructions 
to avoid lost of information at the boundaries of each sector. Also, application of partial backprojections for whole 
object can increase the calculation time. However, partial backprojection algorithms does not need a large oversampling 
to obtain good resolution at a small region of interest and far from the origin. 


\subsection{Backprojection Slice Theorem}

The second algorithm we consider in this paper is based on the Backprojection Slice Theorem. 
Due to the fact that the log-polar coordinates interpolation is not needed for this 
reconstruction, the algorithm is simpler than the above one. Also we note that straight 
usage of the Fourier transform may produce rather big artifacts near the origin (boundary effect), 
caused by the fact that the values on sinogram on the line $s=0$ are not equal to 
$0$. This problem can be solved with usage of short-time Fourier transform:
\begin{equation}
	\hat{f}(t, \sigma) = \int\limits_{-\infty}^{+\infty} f(x) w(x-t) e^{i\sigma x} dx
\end{equation}
with $t=0$. In this work as a function $w$ we are using the Kaiser-Bessel window \cite{kaiser}, 
which can be defined on the mesh $G_{\sf s}$ - see (\ref{tab:symbols}) - using the following formula:
\begin{equation}\label{w}
  w(\beta) = \frac{\Big|I_0\Big(\beta \sqrt{1 - \big(\frac{2i - N_s + 1}{N_s-1} \big)^2} \Big)  \Big|}{|I_0(\beta)|},
\end{equation}
where $I_0(\cdot)$ is a modified Bessel function of the zeroth order. 

In this case, the sequence of steps to retrieve the backprojected image follows:
\begin{description} 
	\item[Step 1.] Transform the source image from sinogram coordinates to semi-polar: $[g]_{\sf s} \to [g]_{\sf P}$.

	\item[Step 2.] For each constant $\theta$:

		\subitem 2.1. Multiply the $s-$axis (polar domain, $s\geq 0)$ of a sinogram with 
		the window function $w(s)$, defined by (\ref{w}): 
		$[\widetilde{g}]_{\sf P}(s,\theta) = [g]_{\sf P}(s,\theta) \cdot w(s) $
		
		\subitem 2.2. Derive \textsc{1d-fft} for the obtained $s-$axis; 

		\subitem 2.3. Multiply the obtained sinogram with the kernel $K^{\sigma} = \frac{1}{\sigma}$ 
		in the frequency domain or its approximation (to avoid division by zero). 

	\item[Step 3.] Interpolate the resulting image to cartesian coordinates, in the frequency domain.
	\item[Step 4.] Apply two-dimensional inverse Fourier transform to obtain the final backprojected image.

\end{description} 

Considering the \textsc{bst} formula, we can notice that the method also has an irregularity at the origin. This irregularity, caused by the division on zero in the 
frequency domain, can be a problem in calculations. The simplest way to avoid this problem is to 
exclude the origin from the calculations. The other way is 
to approximate this division, changing $\sigma = 0$ with some
$\alpha$, where $\alpha>0$ is a small parameter. In this work we use 
\begin{equation}\label{kernel}
	K^{\sigma} = \begin{cases} 
		\frac{1}{\sigma}, & \quad \text{if} \quad  \sigma \neq 0, \\
		\frac{1}{\Delta\sigma}, & \quad \text{otherwise},
	\end{cases}
\end{equation}
where $\Delta\sigma$ is the step of the mesh on $\sigma$.



\section{Regularized FBP: An Application}
\label{sec:reg}

The \textsc{bst} formula (\ref{eq:bst}) can be used to obtain an analytical solution of the 
standard Tikhonov regularization problem in the feature space $U$
\begin{equation}
 \begin{array}{lll} \label{eq:opt}
       \minimize_{f \in U} \|  \Radon f - g \|^2_{L_2} + \lambda \|f\|^2_{L_2} \\
       \end{array}
\end{equation}
In fact, the Euler-Lagrange equations provide the optimality condition for 
the above optimization problem, i.e., $f$ minimizes (\ref{eq:opt}) if and only if \cite{luenberger}
\begin{equation} \label{eq:neq}
(\Radon^* \Radon + \lambda \mathcal I) f (\bm x) = \Radon^* g (\bm x)
\end{equation}
with $\Radon^*$ standing for the adjoint operator of the Radon
transform and $\mathcal I$ the identity operator in $U$. 
In fact, (\ref{eq:neq}) are the so-called normal equations in the 
Hilbert spaces $U$ and $V$. Since $\Radon^* = \Back$ in the usual 
inner-product for $L_2$, the above equation becomes
\begin{equation} \label{eq:neq2}
(\Back \Radon + \lambda \mathcal I) f(\bm x)  = \Back g (\bm x)
\end{equation}
Applying the Fourier transformation on (\ref{eq:neq2}) and using property (\ref{eq:BRf}), 
we obtain the following standard result
\begin{equation} \label{eq:quaseo}
\hat{f}(\bm \omega) \frac{1}{\|\bm \omega\|_2} + \lambda \hat{f}(\bm \omega) = \widehat{\Back g}(\bm \omega) 
\ \ \iff \ \ 
\hat{f}(\bm \omega) \left( \frac{1+\lambda \|\bm \omega\|_2}{ \|\bm \omega\|_2}\right) = \widehat{\Back g}(\bm \omega) 
\end{equation}
From (\ref{eq:quaseo}) is easy to obtain $f$ as a convolution of $\Back g$ with 
an specific two-dimensional filter. If $\lambda=0$ the analytical formula obtained
is exactly the `rho-filter layergram' proposed in \cite{gabor} consisting in a 
post-processing of the backprojection (also mentioned earlier in this manuscript 
as \ \emph{filter of the backprojection}).

The novelty here is that, if we change (\ref{eq:quaseo}) to polar coordinates, we
can immediately apply the \textsc{bst} formula (\ref{eq:bst}). Indeed,
since $\bm \omega = \sigma \bm \xi_\theta$, the pointwise product becomes
\begin{equation}
\hat{f}(\sigma\bm\xi_\theta) \left( \frac{1+\lambda \sigma}{ \sigma}\right) = \widehat{\Back g}(\sigma\bm \xi_\theta) 
\end{equation}
which is essentially the same as
\begin{equation} \label{eq:reg_fst}
\hat{f}(\sigma\bm\xi_\theta) = \left( \frac{1}{1+\lambda \sigma} \right) \hat{g}(\sigma, \theta)
\end{equation}
The above equation is a regularized version of the Fourier-Slice-Theorem and can be 
used to obtain $f$ explicitly through any gridding strategy \cite{marone}.

Applying (\ref{eq:reg_fst}) in the change of variables of the Fourier representation
of $f(\bm x)$ we finally obtain a new representation for the reconstructed image
$f$,
\begin{eqnarray}
f(\bm x) &=& \int_{\R^2} \hat{f}(\bm \omega) e^{i \bm \omega \cdot \bm x} \mathrm d \bm \omega \\
&=& \int_\R \mathrm d\sigma \int_0^{\pi} \mathrm d \theta 
f(\sigma \bm \xi_\theta) |\sigma| e^{i \sigma \bm x \cdot \bm \xi_\theta} \\
&=& \int_\R \mathrm d\sigma \int_0^{\pi} \mathrm d \theta 
\left( \frac{1}{1+\lambda |\sigma|} \right) \hat{g}(\sigma, \theta) |\sigma| e^{i \sigma \bm x \cdot \bm \xi_\theta} \label{eq:reg_fbp}
\end{eqnarray}

Equation (\ref{eq:reg_fbp}) provides exactly the same reconstruction pattern
as a typical filtered backprojection reconstruction algorithm, but with a 
different filter. In fact, we can generalize our regularized strategy
in the following representation
\begin{equation} \label{eq:reg_sol}
f_\lambda(\bm x) = \Back F_\lambda g(\bm x)
\end{equation}
Now, $\{f_\lambda\}$ is a family of solutions of the optimization problem (\ref{eq:opt}),
depending on the regularization parameter $\lambda$. The filter function
$F_\lambda$, in the frequency domain reads
\begin{equation}
\widehat{F_\lambda} (\sigma) = \frac{|\sigma|}{1 + \lambda |\sigma|}
\end{equation}
Our regularized solution (\ref{eq:reg_sol}) depends explicitly on the 
computation of the Backprojection operator $\Back$, and either the \textsc{bst} 
or Andersson's formula can be used.

\section{Numerical Results}
\label{sec:all}

\begin{figure}
	\centering
	\begin{tabular}{cc}
		(a) & (b) 
		\\
		\begin{tikzpicture}[zoomboxarray,
		    zoomboxes below,
		    zoomboxarray columns=1,
		    zoomboxarray rows=1,
		    connect zoomboxes,
		    zoombox paths/.append style={ultra thick, gray}]
		    \node [image node] { \includegraphics[width=0.45\textwidth]{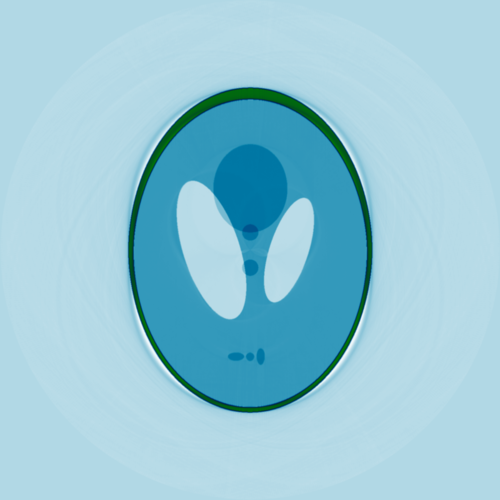} };
		    \zoombox[magnification=3]{0.4,0.4}
		\end{tikzpicture}
		&
		\begin{tikzpicture}[zoomboxarray,
		    zoomboxes below,
		    zoomboxarray columns=1,
		    zoomboxarray rows=1,
		    connect zoomboxes,
		    zoombox paths/.append style={ultra thick, gray}]
		    \node [image node] { \includegraphics[width=0.45\textwidth]{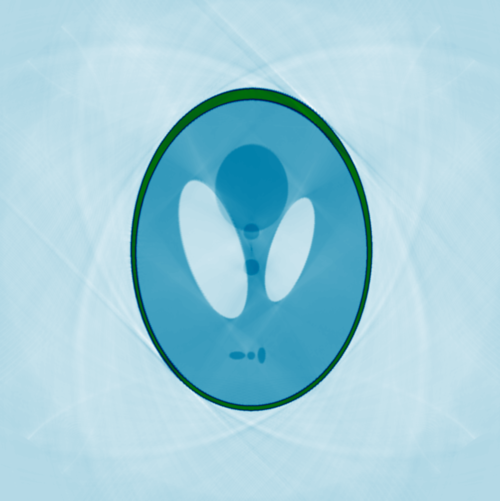} };
		    \zoombox[magnification=3]{0.4,0.4}
		\end{tikzpicture}
	\end{tabular}
	\caption{The results of the log-polar reconstuction using two approaches to cope with origin irregularity: 
	(a),(b) - adaptive $\rho_0$ selection; (c),(d) - partial Backprojections.\label{fig:logpolar_origin}}
\end{figure}

All the algorithms were implemented using the fast Fourier framework \textsc{fftw3} \cite{fftw}. 
We validate our approach using five datasets: two real sinograms and three simulated. The experimental
sinograms (a slice from a wood-fiber and a porous rock) were obtained at the imaging beamline of the 
Brazilian Synchrotron light source and are high-resolution images with $2048\times 1000$ (rays $\times$ angles). 
Therefore, the feature images (either backprojected or filtered-backprojected) were restored with 
$2048\times 2048$ pixels in order
to test the efficiency of the algorithms. The simulated data are: i) the classical shepp-logan phantom
depicted in Figure \ref{fig:phantom_sino}, ii) the circular function of Section \ref{sec:anal} which
has an analytical representation and iii) the following linear combination
\begin{equation}
f(\bm x) = \sum_{j=1}^{1000} \delta (\bm x - \bm a_j)
\end{equation} 
where $\{\bm a_j\}$ are points randomly spanned over the domain $[-0.3,0.3]\times [-0.3,0.3]$,

In section \ref{sec:origin} we described two methods of solving the irregularity near the origin 
for log-polar backprojection. On Fig.\ref{fig:logpolar_origin} we present the comparison of our calculations
using two described methods: on Fig.\ref{fig:logpolar_origin} - the backprojection of the Shepp-Logan test 
function using the adaptive selection of $\rho_0$; on Fig.\ref{fig:logpolar_origin} - log polar reconstruction
with usage of Partial Backprojections. Also we note that on practice the first (adaptive) algoritm works 1.5-2 
times faster. 

In futher tests we compare the results of \textsc{bst} with Log-Polar reconstruction with an adaptive $\rho_0$ selection,  
and the \textsc{nfft} approach \cite{nfft} for the Fourier-Slice-Theorem. 

\begin{figure}
	\begin{tabular}{ccc}
		(a) & (b) & (c) \\
		\begin{tikzpicture} 
	        \begin{axis}[width=0.33\textwidth, compat=newest,
	        ticks=none,enlargelimits=false, axis on top, scale only axis, 
	        axis equal image ] 
	        \addplot graphics 
	        [xmin=0,xmax=256,ymin=0,ymax=256]{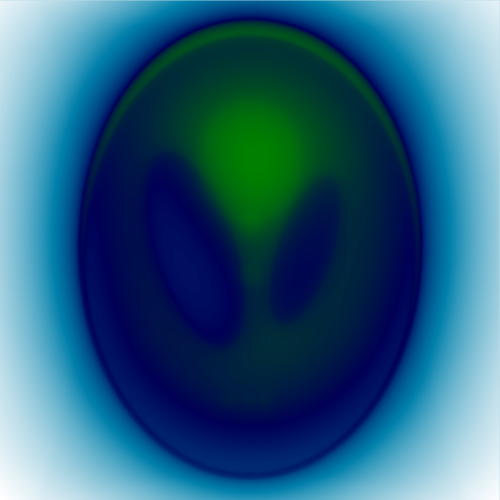};
	        \end{axis}
	 	\end{tikzpicture}	
		&
		\begin{tikzpicture} 
	        \begin{axis}[width=0.33\textwidth, compat=newest,
	        ticks=none,enlargelimits=false, axis on top, scale only axis, 
	        axis equal image] 
	        \addplot graphics 
	        [xmin=0,xmax=256,ymin=0,ymax=256]{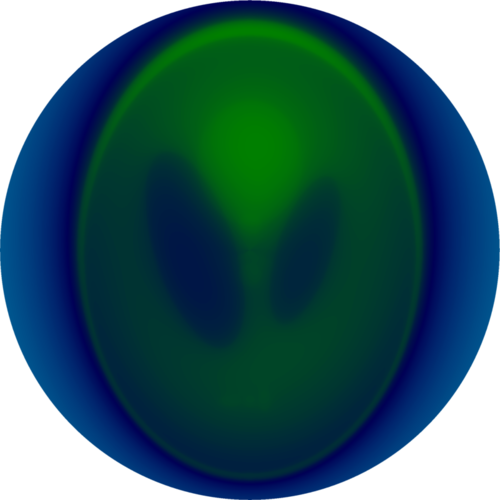};
	        \end{axis}
	 	\end{tikzpicture} 
		&
		\begin{tikzpicture} 
	        \begin{axis}[width=0.33\textwidth, compat=newest,
	        ticks=none,enlargelimits=false, axis on top, scale only axis, 
	        axis equal image] 
	        \addplot graphics 
	        [xmin=0,xmax=256,ymin=0,ymax=256]{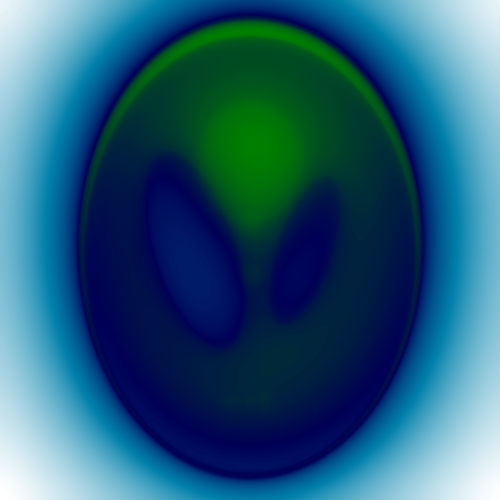};
	        \end{axis}
	 	\end{tikzpicture}
		\\ 
		\begin{tikzpicture} 
	        \begin{axis}[width=0.33\textwidth, compat=newest,
	        ticks=none,enlargelimits=false, axis on top, scale only axis, 
	        axis equal image] 
	        \addplot graphics 
	        [xmin=0,xmax=256,ymin=0,ymax=256]{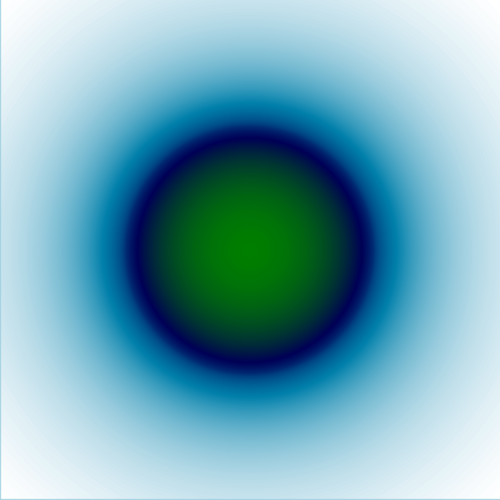};
	        \end{axis}
	 	\end{tikzpicture}
		&
		\begin{tikzpicture} 
	        \begin{axis}[width=0.33\textwidth, compat=newest,
	        ticks=none,enlargelimits=false, axis on top, scale only axis, 
	        axis equal image] 
	        \addplot graphics 
	        [xmin=0,xmax=256,ymin=0,ymax=256]{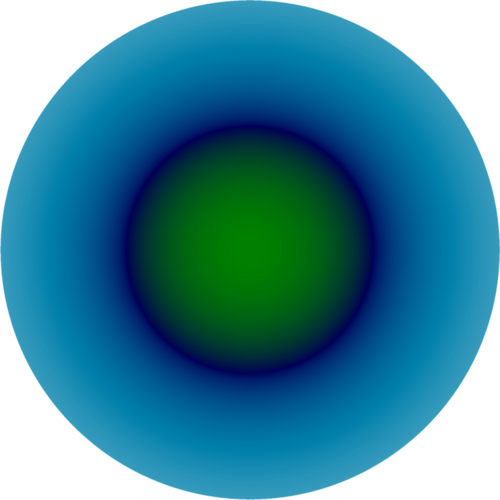};
	        \end{axis}
	 	\end{tikzpicture} 
		&
		\begin{tikzpicture} 
	        \begin{axis}[width=0.33\textwidth, compat=newest,
	        ticks=none,enlargelimits=false, axis on top, scale only axis, 
	        axis equal image] 
	        \addplot graphics 
	        [xmin=0,xmax=256,ymin=0,ymax=256]{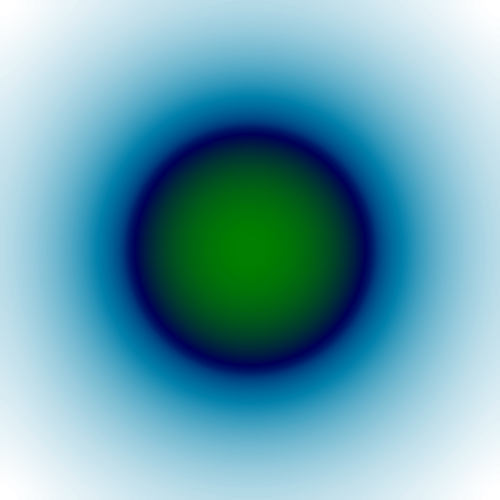};
	        \end{axis}
	 	\end{tikzpicture}
		\\
		\begin{tikzpicture} 
	        \begin{axis}[width=0.33\textwidth, compat=newest,
	        ticks=none,enlargelimits=false, axis on top, scale only axis, 
	        axis equal image] 
	        \addplot graphics 
	        [xmin=0,xmax=256,ymin=0,ymax=256]{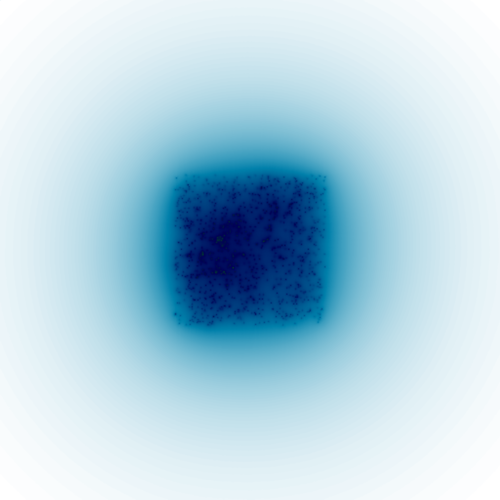};
	        \end{axis}
	 	\end{tikzpicture}
		&
		\begin{tikzpicture} 
	        \begin{axis}[width=0.33\textwidth, compat=newest,
	        ticks=none,enlargelimits=false, axis on top, scale only axis, 
	        axis equal image] 
	        \addplot graphics 
	        [xmin=0,xmax=256,ymin=0,ymax=256]{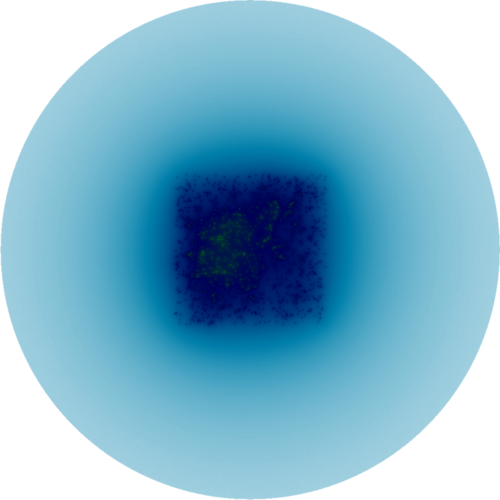};
	        \end{axis}
	 	\end{tikzpicture} 
		&
		\begin{tikzpicture} 
	        \begin{axis}[width=0.33\textwidth, compat=newest,
	        ticks=none,enlargelimits=false, axis on top, scale only axis, 
	        axis equal image] 
	        \addplot graphics 
	        [xmin=0,xmax=256,ymin=0,ymax=256]{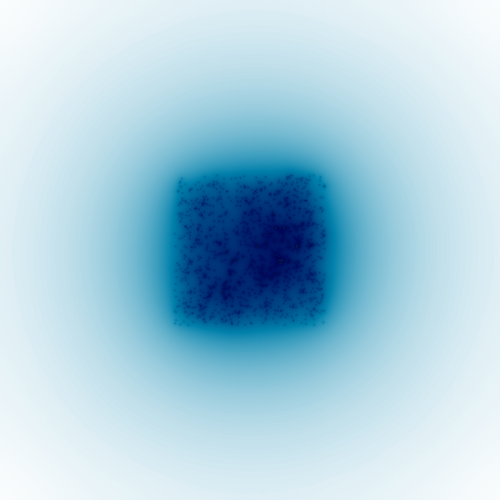};
	        \end{axis}
	 	\end{tikzpicture}
		\\
		\begin{tikzpicture} 
	        \begin{axis}[width=0.33\textwidth, compat=newest,
	        ticks=none,enlargelimits=false, axis on top, scale only axis, 
	        axis equal image] 
	        \addplot graphics 
	        [xmin=0,xmax=256,ymin=0,ymax=256]{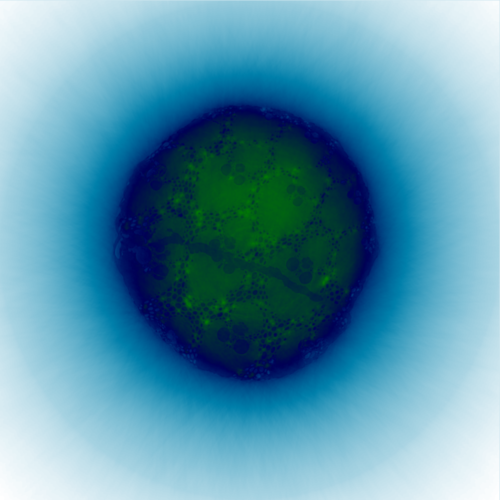};
	        \end{axis}
	 	\end{tikzpicture}
		&
		\begin{tikzpicture} 
	        \begin{axis}[width=0.33\textwidth, compat=newest,
	        ticks=none,enlargelimits=false, axis on top, scale only axis, 
	        axis equal image] 
	        \addplot graphics 
	        [xmin=0,xmax=256,ymin=0,ymax=256]{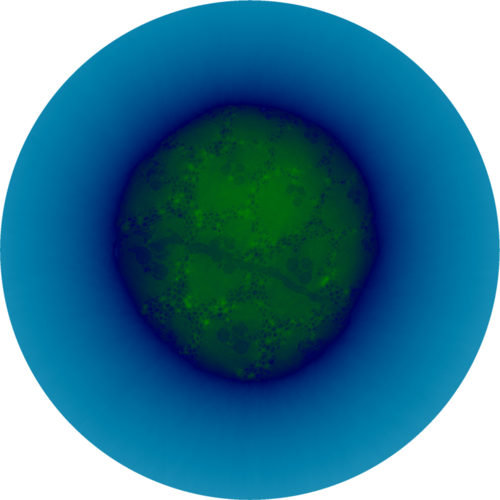};
	        \end{axis}
	 	\end{tikzpicture} 
		&
		\begin{tikzpicture} 
	        \begin{axis}[width=0.33\textwidth, compat=newest,
	        ticks=none,enlargelimits=false, axis on top, scale only axis, 
	        axis equal image] 
	        \addplot graphics 
	        [xmin=0,xmax=256,ymin=0,ymax=256]{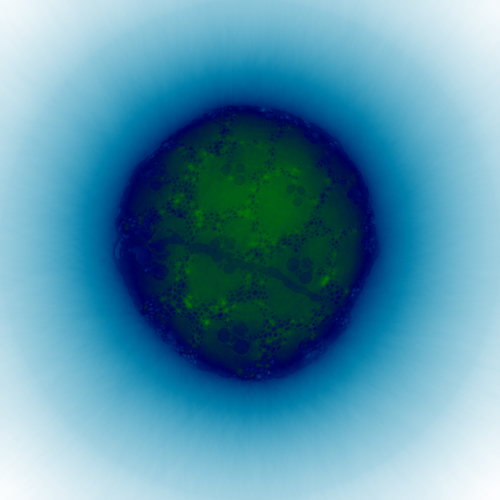};
	        \end{axis}
	 	\end{tikzpicture}
		\\
		\begin{tikzpicture} 
	        \begin{axis}[width=0.33\textwidth, compat=newest,
	        ticks=none,enlargelimits=false, axis on top, scale only axis, 
	        axis equal image] 
	        \addplot graphics 
	        [xmin=0,xmax=256,ymin=0,ymax=256]{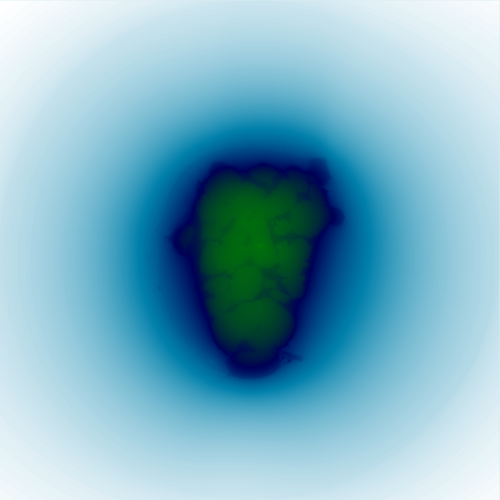};
	        \end{axis}
	 	\end{tikzpicture}
		&
		\begin{tikzpicture} 
	        \begin{axis}[width=0.33\textwidth, compat=newest,
	        ticks=none,enlargelimits=false, axis on top, scale only axis, 
	        axis equal image] 
	        \addplot graphics 
	        [xmin=0,xmax=256,ymin=0,ymax=256]{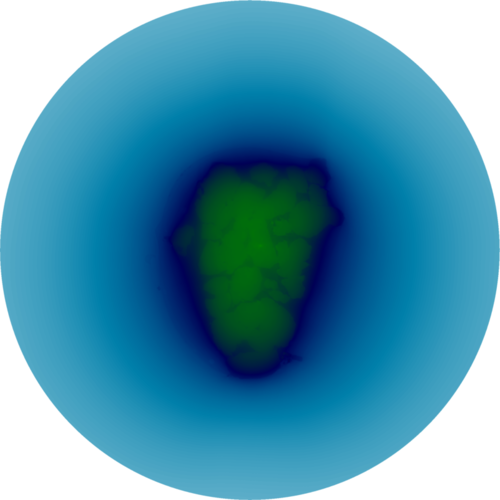};
	        \end{axis}
	 	\end{tikzpicture} 
		&
		\begin{tikzpicture} 
	        \begin{axis}[width=0.33\textwidth, compat=newest,
	        ticks=none,enlargelimits=false, axis on top, scale only axis, 
	        axis equal image] 
	        \addplot graphics 
	        [xmin=0,xmax=256,ymin=0,ymax=256]{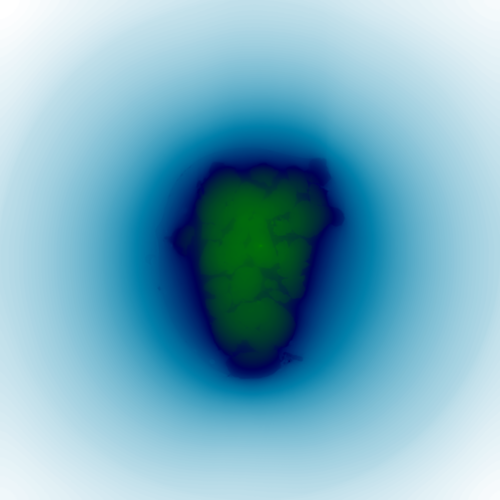};
	        \end{axis}
	 	\end{tikzpicture}
	\end{tabular}
	\caption{\small Comparison between backprojected images. Column (a) shows the results obtained with \textsc{bst}, (b) with Andersson's algorithm and (c) using \textsc{nfft}. See text for details.}
	\label{fig:bp}
\end{figure}

\begin{figure}
	\begin{tabular}{ccc}
		(a) & (b) & (c) \\
		\begin{tikzpicture} 
	        \begin{axis}[width=0.33\textwidth, compat=newest,
	        ticks=none,enlargelimits=false, axis on top, scale only axis, 
	        axis equal image ] 
	        \addplot graphics 
	        [xmin=0,xmax=256,ymin=0,ymax=256]{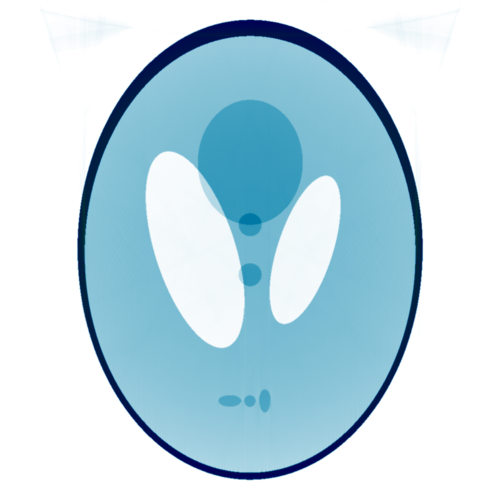};
	        \end{axis}
	 	\end{tikzpicture}	
		&
		\begin{tikzpicture} 
	        \begin{axis}[width=0.33\textwidth, compat=newest,
	        ticks=none,enlargelimits=false, axis on top, scale only axis, 
	        axis equal image] 
	        \addplot graphics 
	        [xmin=0,xmax=256,ymin=0,ymax=256]{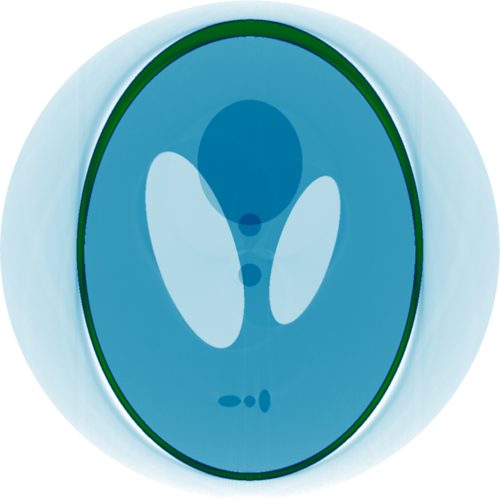};
	        \end{axis}
	 	\end{tikzpicture} 
		&
		\begin{tikzpicture} 
	        \begin{axis}[width=0.33\textwidth, compat=newest,
	        ticks=none,enlargelimits=false, axis on top, scale only axis, 
	        axis equal image] 
	        \addplot graphics 
	        [xmin=0,xmax=256,ymin=0,ymax=256]{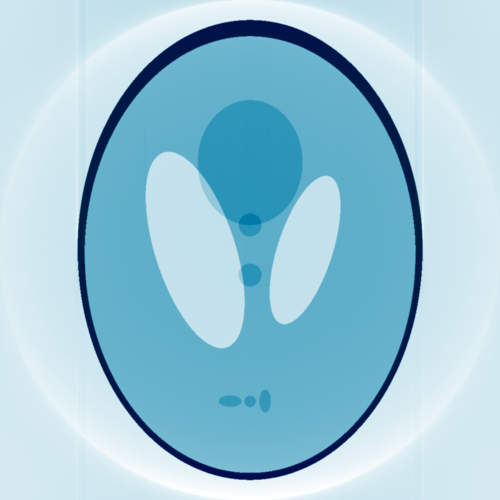};
	        \end{axis}
	 	\end{tikzpicture}
		\\ 
		\begin{tikzpicture} 
	        \begin{axis}[width=0.33\textwidth, compat=newest,
	        ticks=none,enlargelimits=false, axis on top, scale only axis, 
	        axis equal image ] 
	        \addplot graphics 
	        [xmin=0,xmax=256,ymin=0,ymax=256]{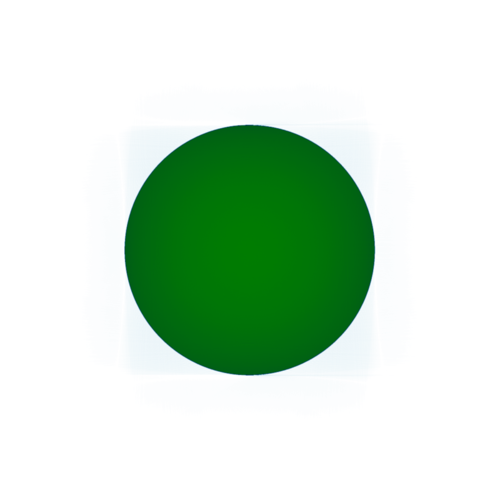};
	        \end{axis}
	 	\end{tikzpicture}	
		&
		\begin{tikzpicture} 
	        \begin{axis}[width=0.33\textwidth, compat=newest,
	        ticks=none,enlargelimits=false, axis on top, scale only axis, 
	        axis equal image] 
	        \addplot graphics 
	        [xmin=0,xmax=256,ymin=0,ymax=256]{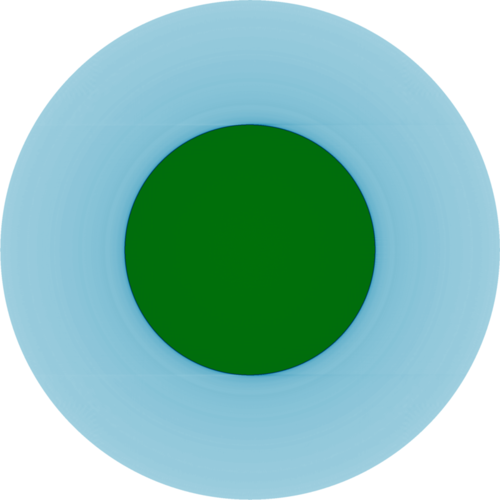};
	        \end{axis}
	 	\end{tikzpicture} 
		&
		\begin{tikzpicture} 
	        \begin{axis}[width=0.33\textwidth, compat=newest,
	        ticks=none,enlargelimits=false, axis on top, scale only axis, 
	        axis equal image] 
	        \addplot graphics 
	        [xmin=0,xmax=256,ymin=0,ymax=256]{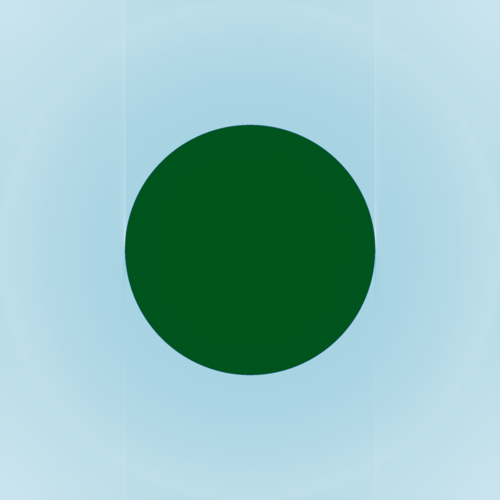};
	        \end{axis}
	 	\end{tikzpicture}
		\\ 
		\begin{tikzpicture} 
	        \begin{axis}[width=0.33\textwidth, compat=newest,
	        ticks=none,enlargelimits=false, axis on top, scale only axis, 
	        axis equal image ] 
	        \addplot graphics 
	        [xmin=0,xmax=256,ymin=0,ymax=256]{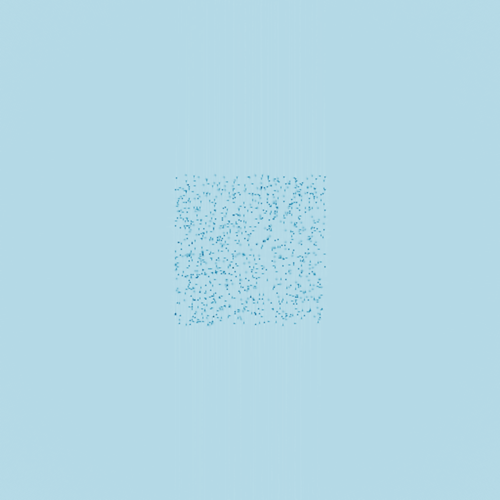};
	        \end{axis}
	 	\end{tikzpicture}	
		&
		\begin{tikzpicture} 
	        \begin{axis}[width=0.33\textwidth, compat=newest,
	        ticks=none,enlargelimits=false, axis on top, scale only axis, 
	        axis equal image] 
	        \addplot graphics 
	        [xmin=0,xmax=256,ymin=0,ymax=256]{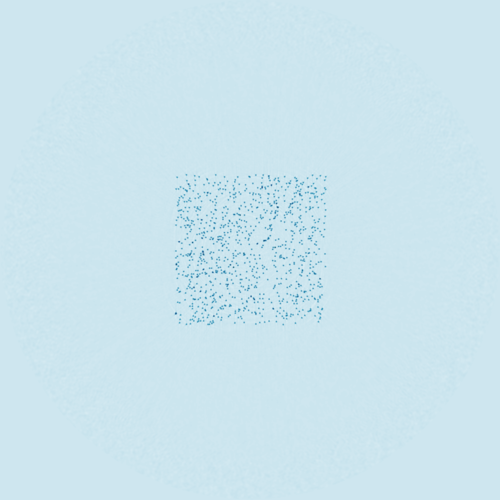};
	        \end{axis}
	 	\end{tikzpicture} 
		&
		\begin{tikzpicture} 
	        \begin{axis}[width=0.33\textwidth, compat=newest,
	        ticks=none,enlargelimits=false, axis on top, scale only axis, 
	        axis equal image] 
	        \addplot graphics 
	        [xmin=0,xmax=256,ymin=0,ymax=256]{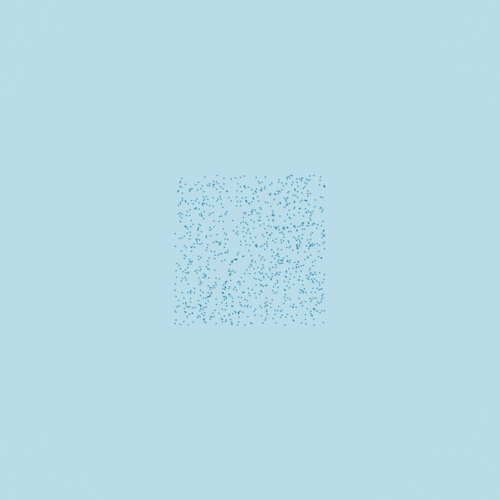};
	        \end{axis}
	 	\end{tikzpicture}
		\\ 
		\begin{tikzpicture} 
	        \begin{axis}[width=0.33\textwidth, compat=newest,
	        ticks=none,enlargelimits=false, axis on top, scale only axis, 
	        axis equal image ] 
	        \addplot graphics 
	        [xmin=0,xmax=256,ymin=0,ymax=256]{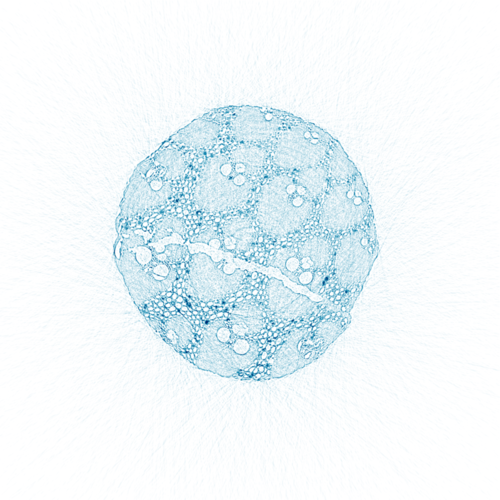};
	        \end{axis}
	 	\end{tikzpicture}	
		&
		\begin{tikzpicture} 
	        \begin{axis}[width=0.33\textwidth, compat=newest,
	        ticks=none,enlargelimits=false, axis on top, scale only axis, 
	        axis equal image] 
	        \addplot graphics 
	        [xmin=0,xmax=256,ymin=0,ymax=256]{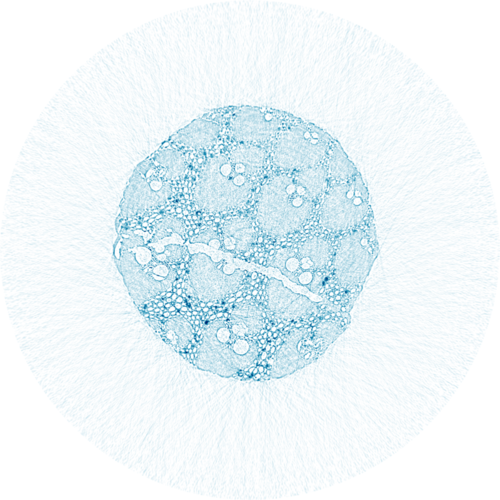};
	        \end{axis}
	 	\end{tikzpicture} 
		&
		\begin{tikzpicture} 
	        \begin{axis}[width=0.33\textwidth, compat=newest,
	        ticks=none,enlargelimits=false, axis on top, scale only axis, 
	        axis equal image] 
	        \addplot graphics 
	        [xmin=0,xmax=256,ymin=0,ymax=256]{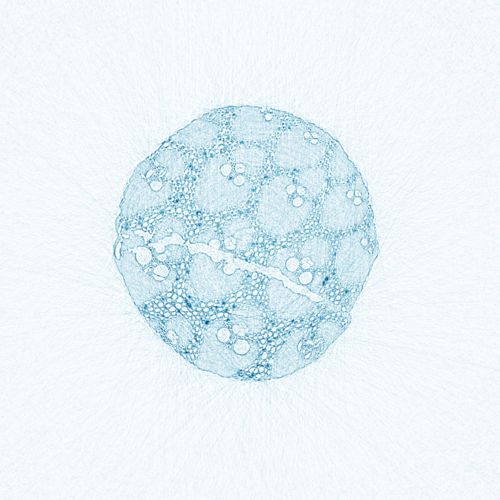};
	        \end{axis}
	 	\end{tikzpicture}
		\\ 
		\begin{tikzpicture} 
	        \begin{axis}[width=0.33\textwidth, compat=newest,
	        ticks=none,enlargelimits=false, axis on top, scale only axis, 
	        axis equal image ] 
	        \addplot graphics 
	        [xmin=0,xmax=256,ymin=0,ymax=256]{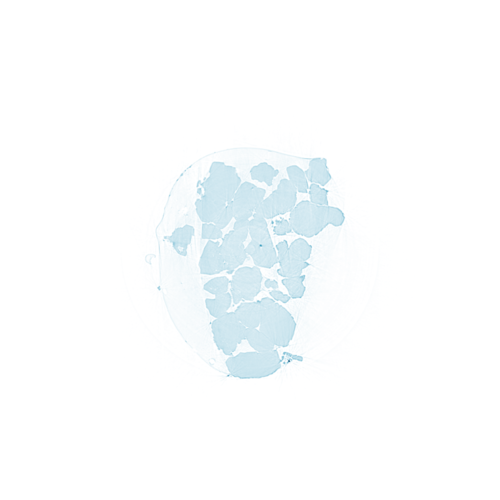};
	        \end{axis}
	 	\end{tikzpicture}	
		&
		\begin{tikzpicture} 
	        \begin{axis}[width=0.33\textwidth, compat=newest,
	        ticks=none,enlargelimits=false, axis on top, scale only axis, 
	        axis equal image] 
	        \addplot graphics 
	        [xmin=0,xmax=256,ymin=0,ymax=256]{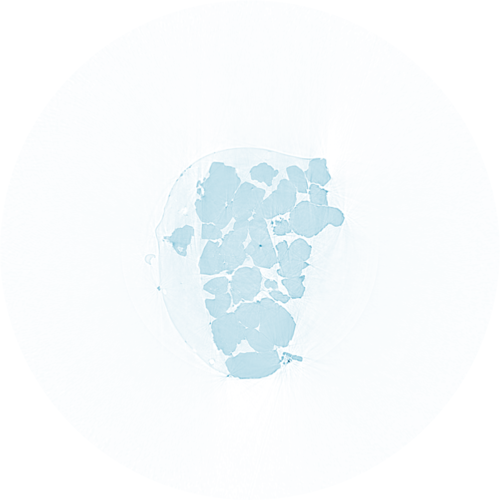};
	        \end{axis}
	 	\end{tikzpicture} 
		&
		\begin{tikzpicture} 
	        \begin{axis}[width=0.33\textwidth, compat=newest,
	        ticks=none,enlargelimits=false, axis on top, scale only axis, 
	        axis equal image] 
	        \addplot graphics 
	        [xmin=0,xmax=256,ymin=0,ymax=256]{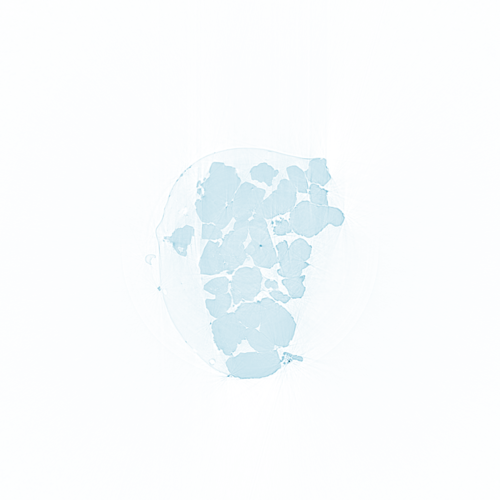};
	        \end{axis}
	 	\end{tikzpicture}
	\end{tabular}
	\caption{Comparison between filtered-backprojected images. Column (a) shows the results obtained with \textsc{bst}, (b) with Andersson's algorithm and (c) using \textsc{nfft}. }
	\label{fig:fbp}
\end{figure}

\begin{figure}
	\centering
	\begin{tabular}{ccc}
	(a) & (b) & (c)
		\\
		\begin{tikzpicture} 
	        \begin{axis}[width=0.33\textwidth, compat=newest,
	        ticks=none,enlargelimits=false, axis on top, scale only axis, 
	        axis equal image ] 
	        \addplot graphics 
	        [xmin=0,xmax=160,ymin=0,ymax=160]{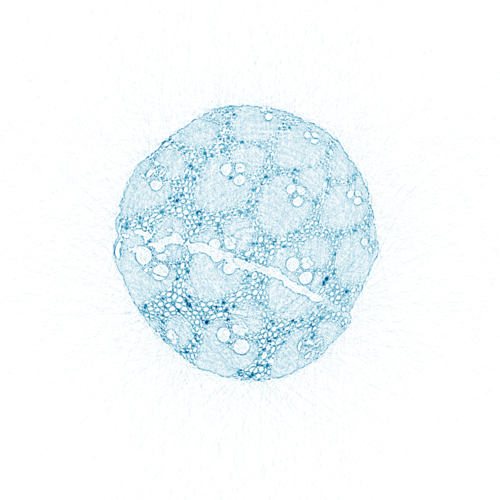};
	        \end{axis}
	 	\end{tikzpicture}
		&
		\begin{tikzpicture} 
	        \begin{axis}[width=0.33\textwidth, compat=newest,
	        ticks=none,enlargelimits=false, axis on top, scale only axis, 
	        axis equal image ] 
	        \addplot graphics 
	        [xmin=0,xmax=160,ymin=0,ymax=160]{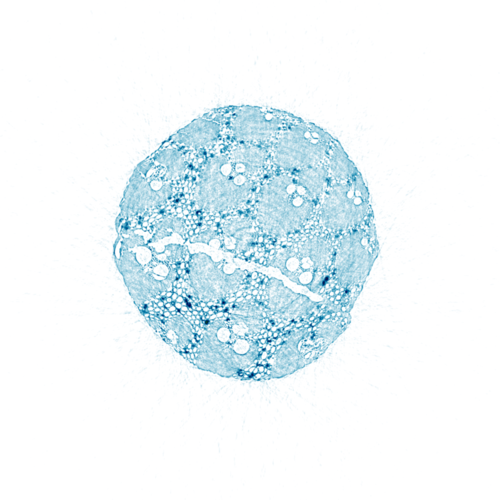};
	        \end{axis}
	 	\end{tikzpicture} 
		&	
		\begin{tikzpicture} 
	        \begin{axis}[width=0.33\textwidth, compat=newest,
	        ticks=none,enlargelimits=false, axis on top, scale only axis, 
	        axis equal image ] 
	        \addplot graphics 
	        [xmin=0,xmax=160,ymin=0,ymax=160]{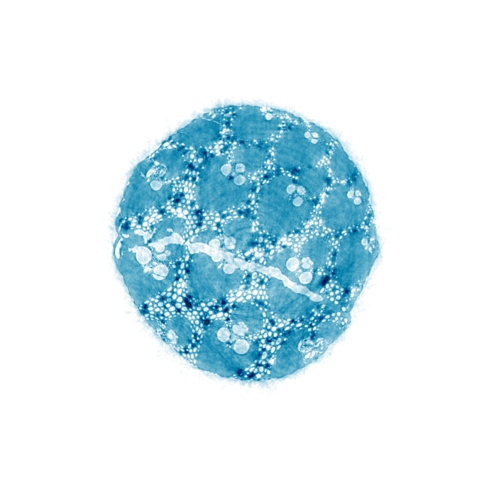};
	        \end{axis}
	 	\end{tikzpicture}		
		\\			
		\begin{tikzpicture} 
	        \begin{axis}[width=0.33\textwidth, compat=newest,
	        ticks=none,enlargelimits=false, axis on top, scale only axis, 
	        axis equal image ] 
	        \addplot graphics 
	        [xmin=0,xmax=160,ymin=0,ymax=160]{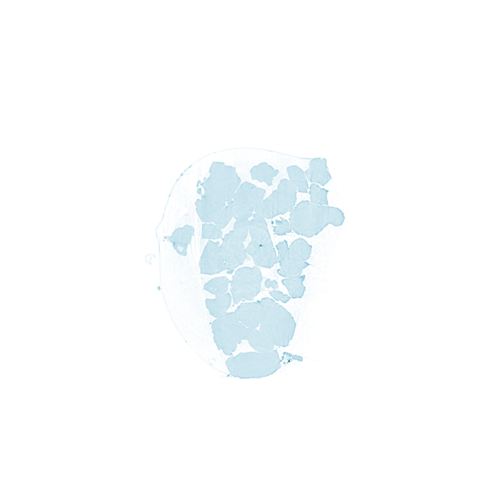};
	        \end{axis}
	 	\end{tikzpicture}
		&
		\begin{tikzpicture} 
	        \begin{axis}[width=0.33\textwidth, compat=newest,
	        ticks=none,enlargelimits=false, axis on top, scale only axis, 
	        axis equal image ] 
	        \addplot graphics 
	        [xmin=0,xmax=160,ymin=0,ymax=160]{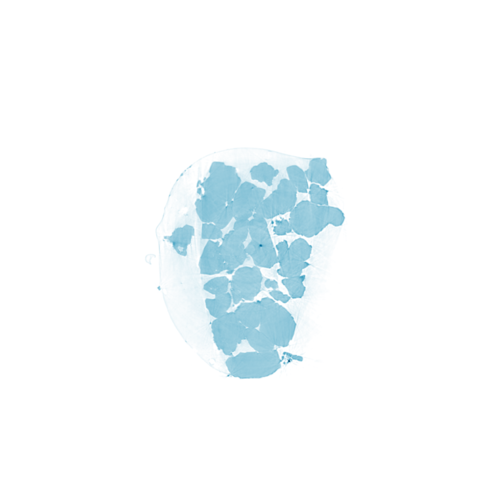};
	        \end{axis}
	 	\end{tikzpicture} 
		&	
		\begin{tikzpicture} 
	        \begin{axis}[width=0.33\textwidth, compat=newest,
	        ticks=none,enlargelimits=false, axis on top, scale only axis, 
	        axis equal image ] 
	        \addplot graphics 
	        [xmin=0,xmax=160,ymin=0,ymax=160]{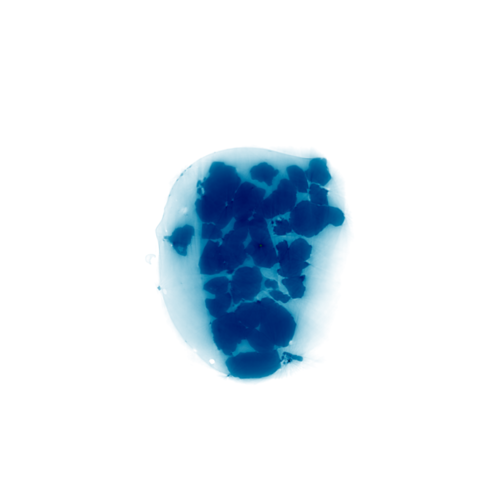};
	        \end{axis}
	 	\end{tikzpicture} 
	\end{tabular}
	\caption{Reconstruction with real data and the regularized filtered-backprojection 
described in Section \ref{sec:reg} for different $\lambda's$ and using \textsc{bst}. 
(a) $\lambda=0.002$, (b) $\lambda=0.02$ and (c) $\lambda=0.2$.}
	\label{fig:reg}
\end{figure}

The regularized filtered-backprojection algorithm described in Section \ref{sec:reg} was applied to the noisy
data gathered for the wood-fiber and the rock sample. The results are shown 
in Figure \ref{fig:reg} for the wood-fiber and the rock sample using only three values for the 
regularization parameter $\lambda$. In fact, an algorithm for the selection of the optimal parameter is beyond the scope
of this manuscript. Regulared filtered-backprojected images were obtained using $\lambda \in \{0.002,0.02,0.2\}$.
As it is known from Tikhonov regularization schemes, the bigger is $\lambda$, the smoother the resulting
image will be. This is clearly visible in Figure \ref{fig:reg}, what indicates that such an approach
could be used to increase the constrast in the final reconstructed image.

All algorithms are fast due to usage of convolutions and Fast Fourier techniques. Computational 
complexity - that we denote by $\Omega$ - of Andersson's approach is similar to the complexity of the two dimensional 
convolution, i.e. 
\begin{equation}
	\Omega_{Andersson} = N_{\theta} N_{\rho} (\log_2 N_{\rho} + \log_2 N_{\theta})+ \Omega_{\sf L},
\end{equation}
where $\Omega_{\sf L}$ is the summarized complexity of all log-polar interpolations (sinogram to log-polar and
log-polar to cartesian). Here we assume that the Fourier transform of the kernel was pre-calculated 
(numerically, or analytically, like in \cite{andersson}) and is not taking it into account for $\Omega_{Andersson}$.
In case of adaptive $\rho_0$ selection, $N_{\rho}$ is being calculated using (\ref{N_rho}), which leads
to some oversampling. For most common sizes of the sinogram (i.e. $512 < N_s, N_x, N_y < 16000$ ) this oversampling 
is not very big ($N_{\rho} \approx 4N_s$), and the complexity of log-polar reconstruction can be estimated as 
$8N^2(\log_2N+1)$.

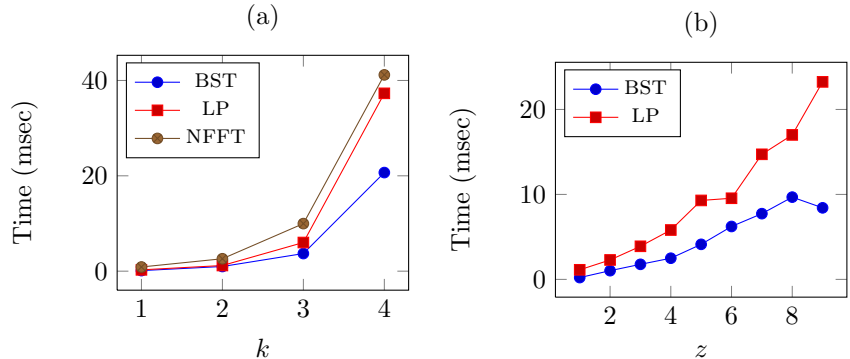
\begin{figure}
	\centering
	\begin{tabular}{ccc}
		\begin{tikzpicture}
		\begin{axis}[width=0.45\textwidth,
		  title={(a)},
		  xlabel=$k$,
		  ylabel=Time (msec),
		  legend pos=north west]
		  \addplot table [x index=0, y index=1]{times.dat};
		  \addlegendentry{\textsc{bst}}
		  \addplot table [x index=0, y index=2]{times.dat};
		  \addlegendentry{\textsc{lp}}
		  \addplot table [x index=0, y index=3]{times.dat}; 
		  \addlegendentry{\textsc{nfft}}
		\end{axis}
		\end{tikzpicture}
		&
		\begin{tikzpicture}
		\begin{axis}[
		  title={(b)},
		  width=0.45\textwidth,
		  xlabel=$z$,
		  ylabel=Time (msec),
		  legend pos=north west]
		  \addplot table [x index=0, y index=1]{zp.dat};
		  \addlegendentry{\textsc{bst}}
		  \addplot table [x index=0, y index=2]{zp.dat};
		  \addlegendentry{\textsc{lp}}
		\end{axis}
		\end{tikzpicture} 		
	\end{tabular}
	\caption{(a) Comparison of reconstruction time $\tau$ {\it versus} size of the image $N=256\times 2^k$
	with a constant zero-padding. (b) Reconstruction time {\it versus} zero-padding coefficient
	$z$, transforming the domain $s\in [0,1]$ to $s\in [0,z]$ in polar coordinates.}
	\label{fig:performance}
\end{figure}

The complexity of reconstruction, based on \textsc{bst} formula also depends on the size of the final image. 
Much of computational complexity falls on the one-dimensional Fourier transforms, whose size is equal to 
the number of rays in the sinogram, and to final 2D Fourier transform in cartesian coordinates.  
Hence, the complexity is 
\begin{equation}
	\Omega_{{\sf BST}} = N_{\theta}N_s \log_2(N_s) + N_xN_y(\log_2N_x + \log_2N_y)+ \Omega_{\sf P}, 
\end{equation}
where $\Omega_{\sf P}$ denotes the total amount of complexity for polar 
interpolations (polar to cartesian), usually $\Omega_{\sf P} = O(N_s^2)$. Note
that for clear reconstruction we also need some oversampling on $s$, but 
not so big as in the previous case (not more than two times). The comparison between time for
reconstruction, depending on the size of the sinogram (for model task of Shepp-Logan phantom reconstruction) 
is presented on Figure \ref{fig:performance}.(a). On Figure \ref{fig:performance}.(b) we present the dependence of the reconstruction time \emph{versus} the zero-padding factor, say $z$, for the log-polar (adaptive) and \textsc{bst}. 
In our simulations, the zero-padding $z$ increase the support of the sinogram from $[0,1]$ to 
$[0,z]$ in polar coordinates. 

The backprojection (or reconstruction) for a slice with size $1024 \times 2048$ (polar coordinates, rays $\times$ angles) 
can be obtained in about $690$ milliseconds on a modest computer using \textsc{cpu} Intel(R) Core(TM) i7-3770 CPU @ 3.40GHz
with only $8$ threads. Of course, the programs developed by authors can be sufficiently optimized and 
powered for \textsc{gpu}, which can make the backprojection considerably faster using either 
Andersson's formula or \textsc{bst}.

\begin{figure}
	\centering
	\begin{tikzpicture}
		\begin{axis}[width=0.8\textwidth,
		  title={(a)},
		  xlabel=Sinogram error (\%),
		  ylabel=Filtered backprojection error (\%),
		  legend pos=north west]
		  \addplot table [x index=0, y index=1]{mse.dat};
		  \addlegendentry{\textsc{bst}}
		  \addplot table [x index=0, y index=2]{mse.dat};
		  \addlegendentry{\textsc{lp}}
		  \addplot table [x index=0, y index=3]{mse.dat}; 
		  \addlegendentry{\textsc{nfft}}
		  \addplot table [x index=0, y index=4]{mse.dat}; 
		  \addlegendentry{\textsc{slan}}
		  \addplot table [x index=0, y index=5]{mse.dat}; 
		  \addlegendentry{\textsc{bres}}
		\end{axis}
	\end{tikzpicture}
	\caption{(a) Mean square error (MSE) of the result, obtained using different algorithms, in dependence 
	on MSE of enter data (sinogram).}
	\label{fig:accuracy}
\end{figure}
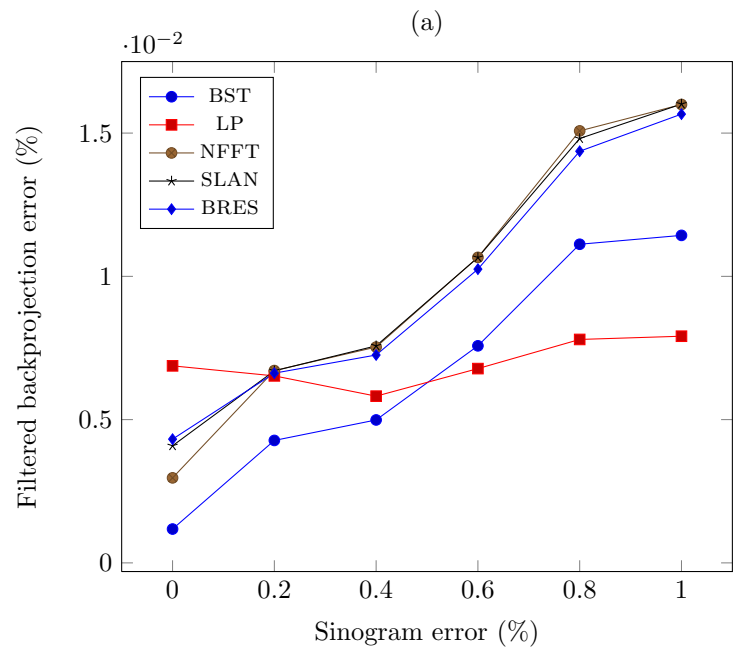

Figure \ref{fig:accuracy} presents some benchmarks of accuracy for the developed algorithms with other known backprojection techniques. 
More precisely, we present the resulting mean squared error (MSE) \emph{versus} the error in the input data (i.e., the sinogram). 
Calculation were done for the Shepp-Logan phantom with addition of Poisson noise to the analytical sinogram. From 
Figure \ref{fig:accuracy} one can note that for weak noises BST shows the best accuracies, while Log-Polar reconstruction is 
very stable to strong noise.

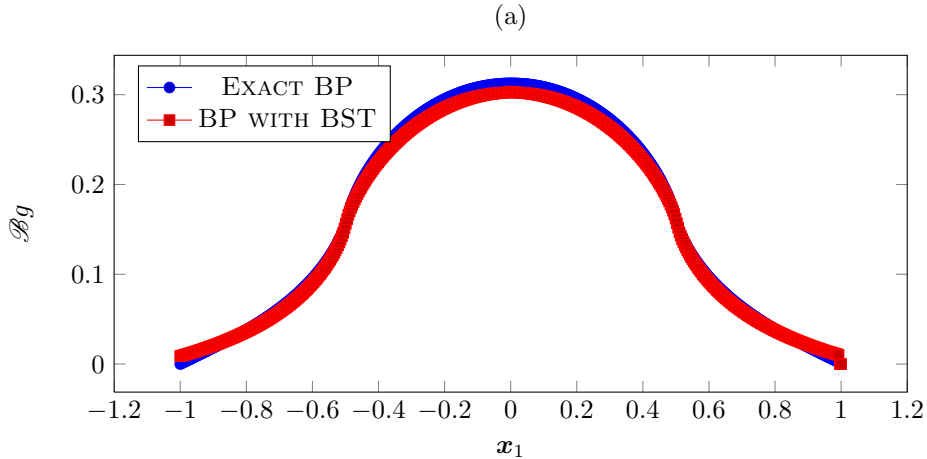
\begin{figure}
	\centering
	\begin{tikzpicture}
		\begin{axis}[width=1.0\textwidth,height=0.5\textwidth,
		  title={(a)},
		  xlabel= $\bm x_1$,
		  ylabel= $\Back g$,
		  legend pos=north west]
		  \addplot table [mark=none, x index=0, y index=1]{slice2.dat};
		  \addlegendentry{\textsc{Exact BP}}
		  \addplot table [mark=none, x index=0, y index=2]{slice2.dat};
		  \addlegendentry{\textsc{BP with BST}}
		\end{axis}
	\end{tikzpicture}
	\caption{Analytical and numerical (BST) backprojections for Example 2, a circular function 
	presented in Section \ref{sec:anal}}.
	\label{fig:slice}
\end{figure}

On Figure \ref{fig:slice} we present the slice of the reconstruction of our analytical example (see Section \ref{sec:anal}, Example 2). 
Since we obtain backprojection of the circle function analytically, we can compare the result of numerical backprojection with BST 
and analytical solution.

\section{Conclusion}
\label{sec:conc}

In this manuscript we have proposed a new backprojection technique (\textsc{bst}) and compared it against two
other already established algorithms, the log-polar (\textsc{lp}) approach from Andersson \cite{andersson} and 
with the use of Non-Uniform Fast Fourier Transforms (\textsc{nfft}) from \cite{nfft}. With the increasing
size of input data, measured at imaging beamlines from synchrotron facilities, the need for fast 
post-processing of the data becomes an immediate demand. Either conventional reconstruction
techniques like filtered backprojection or more robust iterative methods can benefit from a 
fast backprojection algorithm. Finally, in order to demonstrate this possibility in practice, we have provided an 
application where the Tikhonov image is computed by means of the new \textsc{bst} formula in a short amount of time, 
thereby enabling the practical application of interesting regularization schemes to very large datasets.

\appendix

\section{Integral representations}
\label{app:int}

We use the following standard representation for path integrals, the proof can be found in \cite{helgason}: \emph{for a continuously differentiable $m \colon \R^m \to \R$
such that $\|\nabla m \| \not= \bm{0}$ it is true that}:
\begin{equation}\label{eq:pathint}
\int_{\R^m} h(\bm{y}) \delta(m(\bm{y})) \mathrm d\bm{y} =  \int_{C^{-1}(0)} \frac{h(\bm{y}(s))}{\|\nabla m(\bm{y}(s))\|} \mathrm ds (\bm{x}).
\end{equation}
where $\mathrm ds(\bm x)$ is the arclength measure along curve $C^{-1}(0)$. Assuming that $g\in V$ is a given sinogram and $\bm{x}$ a pixel in the reconstruction region. 
The backprojection (\ref{adjoint}) of $g$ is defined as the contribution of all 
possible straight lines, parameterized by the angle $\theta$, and passing through $\bm{x}$.
Using the sifting property of the Delta distribution, we have
\begin{equation}
\Back g(\bm{x}) = \int_0^{\pi} g(\bm{x}\cdot \bm{\xi}_\theta,\theta) \mathrm d \theta 
= \int_0^{\pi} \int_\R g(t,\theta) \delta (t-\bm{x}\cdot\bm{\xi}_\theta) \mathrm d t \mathrm d \theta 
\end{equation}
Switching the above integral from $(t,\theta)$ coordinates to cartesian coordinates $\bm{y}\in\R^2$ we
have $|t|\mathrm d t \mathrm d \theta = \mathrm d\bm{y}$; where $\Back g$ now becomes 
\begin{equation} \label{eq:shit}
\Back g(\bm{x}) = \int_{\R^2} [g]_{\mathsf c}(\bm{y}) \delta ( m(\bm{y}) ) \frac{1}{\|\bm{y}\|}\mathrm d\bm{y}
\end{equation}
with $[g]_{\mathsf c}(\bm{y}) = g(t(\bm{y}), \theta(\bm{y}))$ refering to the sinogram $g$
in cartesian coordinates. In fact, $|t| = \|\bm{y}\|$ is the unsigned distance to the origin
and $\theta = \arctan(\frac{\bm{y}_2}{\bm{y}_1}) \in [0,\pi]$ is 
the angle with respect to the $\bm{y}_1$-axis. Function $m$ reads
\begin{eqnarray}
m(\bm{y}) &=& t - \bm{x}\cdot \bm{\xi}_\theta = \|\bm{y}\|- \bm{x}_1\cos \theta(\bm y)- 
\bm{x}_2 \sin  \theta(\bm y)  \\
&=& \|\bm{y}\| - \bm{x}_1 \frac{\bm y_1}{\|\bm y\|} - \bm{x}_2 \frac{\bm {y}_2}{\|\bm y\|} 
= \|\bm{y}\| - \frac{(\bm{x}_1\bm{y}_1 + \bm{x}_2 \bm{y}_2) }{\|\bm{y}\|}\\
&=& \frac{\bm{y} \cdot (\bm{y} - \bm{x})}{\|\bm{y}\|} \label{eq:sigma}
\end{eqnarray}
From (\ref{eq:sigma}), (\ref{eq:shit}) and the property $\delta(au) = \frac{1}{|a|}\delta(u)$ for all $u \in \R$, 
the backprojection now follows:
\begin{equation} \label{eq:backCirc}
\Back g(\bm{x}) = \int_{\R^2} [g]_{\mathsf c}(\bm{y}) \delta( \kappa_{\bm{x}}(\bm{y})) \mathrm d\bm{y}, \ \ \ \ 
\kappa_{\bm{x}}(\bm{y}) = \bm{y} \cdot (\bm{y} - \bm{x})
\end{equation}
It should be noted that, for a fixed $\bm{x}\in\R^2$, the set $\kappa^{-1}_{\bm{x}}(0) = \{\bm{y}\in\R^2: \kappa_{\bm{x}} (\bm{y}) = 0\}$ 
is defined as a circle in the plane. Indeed, since
$\bm{y} \cdot (\bm{y} - \bm{x}) = \bm{y}\cdot\bm{y} - 2 \bm{y}\cdot\left(\frac{\bm{x}}{2}\right) = 
\left\|\bm{y}-\frac{\bm{x}}{2}\right\|^2 - \left\|\frac{\bm{x}}{2}\right\|^2$,
it follows that $\kappa^{-1}_{\bm{x}}(0)$ is a circle passing through the origin $\bm{y}=\bm{0}$,
centered at $\frac{1}{2}\bm{x}$ and with radius $\frac{1}{2}\|\bm{x}\|$. Since $\kappa_{\bm{x}}^{-1}(0) = \{\frac{1}{2}\bm{x} + r\bm{\xi}_\theta: \theta \in [0,2\pi], \ r=\frac{1}{2}\|\bm{x}\|\}$ is a parametric representation of the circle, the backprojection operator also reads, in an alternative form: {\it $\Back$ is a stacking operator through circles $\kappa_{\bm x}^{-1}(0)$:}
\begin{equation}
\Back g(\bm{x}) = \int_{\kappa^{-1}_{\bm{x}}(0)} \frac{[g]_{\mathsf c}(\bm{y})}{\|2\bm{y}-\bm{x}\|} \mathrm ds 
= \frac{1}{2} \int_{0}^{2\pi} [g]_{\mathsf c}\left(\frac{1}{2}\bm{x} + \frac{1}{2}\|\bm{x}\| \bm{\xi}_\theta \right) 
\mathrm d\theta
\end{equation}
The above representation follows from $ds = \frac{1}{2}\|x\|d\theta$, (\ref{eq:backCirc}) and (\ref{eq:pathint}) with  $\nabla \kappa_{\bm{x}}(\bm{y})=2\bm{y}-\bm{x}$.
Last equality comes from $\bm{y} = \frac{1}{2}\bm{x} + \frac{1}{2}\|\bm{x}\|\xi_\theta \in \kappa^{-1}_{\bm{x}}(0)$ 
for some $\theta$. Therefore, in cartesian coordinates, the backprojection contribution for a ball 
$\{\bm{z}\in\R^2: \|\bm{z}-\bm{x}\|\leq \epsilon \}$ 
comes from a family of circles passing through the ball and the origin, this fact seems to be related to the \emph{comet-tail region} mentioned by \cite{faridani}. 

\section{Point Source in Log-Polar Coordinates}
\label{sec:psf}

In this section we present the details for the following log-polar representation of the backprojected image
\begin{eqnarray}
b(e^\rho \bm \xi_\theta) &=& \int \mathrm d \beta \int \mathrm d u \ \underbrace{\delta( \bm a\cdot\bm \xi_\beta - e^u}_{m(u)} ) \delta(\underbrace{ 1 - e^{\rho-u} \cos (\theta-\beta)}_{\ell(u)}) \\
&=& \frac{1}{\sqrt{(\cos(\theta-\phi) e^\rho)^2 + (e^A - \sin(\theta-\phi) e^\rho)^2}}
\end{eqnarray}

\begin{proof}
Using the property of the Delta distribution \cite{bracewell},
\begin{equation} \label{eq:delta}
\delta( \ell(u)) = \sum_{k} \frac{\delta (u - u_k)}{|\ell'(u_k)|}
\end{equation}
where ${u_k}$ are roots of $\ell(u)$. In our case, since $\ell$ has only one zero
\begin{equation}
u_0 = \rho + \ln \cos(\theta-\beta)
\end{equation}
and $\ell'(u_0) = 1$, function $b$ reads
\begin{equation}
b(e^\rho \bm \xi_\theta) = \int \mathrm d \beta \ m(u) \delta(u-u_0) 
\end{equation}
Due to the sifting property of the Delta, we obtain
\begin{eqnarray}
b(e^\rho \bm \xi_\theta) &=&  \int \mathrm d \beta \ m( \rho + \ln \cos(\theta-\beta) ) \\
&=& \int \mathrm d \beta \ \delta( \bm a \cdot \bm \xi_\beta - e^{\rho + \ln \cos(\theta-\beta)} ) \\
&=& \int \mathrm d \beta \ \delta( \bm a \cdot \bm \xi_\beta - e^\rho \cos(\theta-\beta) ) \\
\end{eqnarray}
Using $\bm a = e^A \bm \xi_\phi$ we obtain
\begin{eqnarray}
b(e^\rho \bm \xi_\theta) &=& \int \mathrm d \beta \ \delta( e^A \bm \xi_\phi \cdot \bm \xi_\beta - e^\rho \cos(\theta-\beta) ) \\
&=& \int \mathrm d \beta \ \delta( \underbrace{e^A \cos(\phi-\beta) -  e^\rho \cos(\theta-\beta)}_{L(\beta)} ) \label{eq:quase}
\end{eqnarray}
Since $\cos(\theta-\beta) = \cos((\phi-\beta) + (\theta-\phi))$,  function $L$ can be rewritten as
\begin{eqnarray*}
L(\beta) &=& e^A \cos(\underbrace{\phi-\beta}_{\alpha}) - e^\rho\left[\cos(\phi-\beta) \underbrace{\sin (\theta-\phi)}_{S} - \sin (\phi-\beta) \underbrace{\cos (\theta-\beta)}_{C}  \right] \\
&=& e^A \cos \alpha - e^\rho S \cos \alpha + e^\rho C \sin \alpha \\
&=& \underbrace{(e^A - S e^\rho)}_{z} \cos \alpha + \underbrace{e^\rho C}_{y} \sin \alpha
\end{eqnarray*}
Hence, the root $\beta_0$ of $L$ must satisfy $L(\beta_0) = 0$, i.e., $\tan \alpha = - \frac{z}{y} \equiv k$
and $\sin \alpha = \frac{k^2}{1+k^2}$, $\cos \alpha = \frac{1}{1+k^2}$. Therefore,
\begin{eqnarray}
L'(\beta_0) &=& z \frac{k}{\sqrt{1+k^2}} - y \frac{1}{\sqrt{1+k^2}} 
= - \frac{z^2}{\sqrt{y^2 + z^2}} - \frac{k^2}{\sqrt{y^2 + z^2}} \\
&=& - \sqrt{y^2 + z^2} = - \sqrt{(C e^\rho)^2 + (e^A - S e^\rho)^2} \label{eq:quase2}
\end{eqnarray}
Finally,
\begin{eqnarray} 
b(e^\rho \bm \xi_\theta) &=& \int \mathrm d \beta \frac{\delta( \beta - \beta_0) )}{|L'(\beta_0)|} = \frac{1}{|L'(\beta_0)|} \\
&=& \frac{1}{\sqrt{(C e^\rho)^2 + (e^A - S e^\rho)^2}} \\
&=& \frac{1}{\sqrt{(\cos(\theta-\phi) e^\rho)^2 + (e^A - \sin(\theta-\phi) e^\rho)^2}}
\end{eqnarray}
which is the final representation of the backprojected image in log-polar coordinates.
\end{proof}


\bigskip

\bibliography{main.bib}{}
\bibliographystyle{unsrt}

\end{document}